\documentclass[leqno,a4paper,titlepage,twoside,10pt]{smfart}

\usepackage[T1]{fontenc}

\usepackage[francais]{babel}

\usepackage{bbm,eucal,mathrsfs,graphicx,epsfig,latexsym,eufrak,pifont,multirow,color,rotating}

\usepackage{smfthm,bull}

\usepackage[all]{xy}

\newcounter{infra}[page]

\newenvironment{dem}[1][]{%
{\bf D\'emonstration #1 : }}{%
\hspace*{\fill}\nolinebreak[1]\hspace*{\fill}\underline{\bf Q.e.d.}\\}

\newenvironment{dem*}[1][]{%
{\bf D\'emonstration #1 : }}{%
 }

\newenvironment{eq*}{\begin{eqnarray*}}{\end{eqnarray*}}

\newenvironment{liste}{\begin{itemize}}{\end{itemize}} 

\newenvironment{rem}{{\bf{Remarque : }}}{}

\newtheorem{defin}{D\'efinition}

\newtheorem{thm}{Th\'eor\`eme}[section] \newtheorem{lem}[thm]{Lemme}
\newtheorem{pro}[thm]{Proposition}

\newcommand{\croi}{\times}

\newcommand{\adh}{\overline}

\newcommand{\antidiago}{.\hskip -.5ex \cdot \hskip -.5 ex \raisebox{1 ex}.}

 \newcommand{\boufl}{\xymatrix{ \!\!
\ar@(ur,dr) }}

\newcommand{\cad}{\mbox{\it c-\`a-d }}
\newcommand{\call}{\mathscr}

\newcommand{\cd}{\rangle}
\newcommand{\cf}{{\it cf. }}

\newcommand{\cg}{\langle}

\newcommand{\ch}{\vee} 
 
 \newcommand{\codim}{\mathrm{codim}}

 \newcommand{\con}{\supseteq}

\newcommand{\diag}{\mathrm{diag}} 
\newcommand{\dij}{\sqcup} 
\newcommand{\Dij}{\bigsqcup}

\newcommand{\axaa}{\sta{\alpha}{\bullet}\;\;\sta{\alpha'}{\bullet}}

\newcommand{\spii}{\sta{\alpha_1}{\circ} \!-\! \sta{\alpha_2}{\circ} \!\Leftarrow \!\sta{\alpha_3}{\bullet}}

\newcommand{\dnn}{\sta{\alpha_1}{\bullet} \!-\! \sta{\alpha_2}{\circ} \!- ... -\!\sta{\alpha_{n-2}}{\circ}\!  <\! \raisebox{-1 ex}{$\circ_{\alpha_{n-1}}$}\hskip -5.5ex\raisebox{1 ex}{$\circ^{\alpha_n}$}}
\newcommand{\donne}{\mapsto}

 \newcommand{\ens}{\mathbbm}

\newcommand{\equi}{\Leftrightarrow} \newcommand{\et}{\mbox{ et }}

\newcommand{\exist}{\exists\:}

 \newcommand{\gamb}{\leadsto}

   \newcommand{\goth}{\mathfrak}
\newcommand{\Hom}{\mathrm{Hom}}

 \newcommand{\ie}{\emph{i.e. }}
 \newcommand{\impliq}{\Rightarrow}

 \newcommand{\infi}{\infty}
\newcommand{\infra}[1]{\stepcounter{infra}\setcounter{footnote}{\value{infra}}%
\footnote{#1}}

\renewcommand{\int}{\mathrm{int\:}}

\newcommand{\inter}{\bigcap} \newcommand{\inv}{^{-1}}
 \newcommand{\iso}{\simeq}

\newcommand{\kk}{{\mathbf{k}}}

\newcommand{\limi}{\ma \lim_{a \to\infi}}

\newcommand{\limz}{\ma \lim_{ a \to 0}}

\newcommand{\ma}{\displaystyle} 

\newcommand{\moins}{\:\setminus\:} 

 \newcommand{\Ou}{\mbox{\Ou}}
\newcommand{\ou}{\mbox{ ou }} 
\newcommand{\p}{\:\:.}  
\newcommand{\PP}{\mathbbm{P}}

\newcommand{\Pic}{\mathrm{Pic}} 
\newcommand{\pic}{\mathrm{pic}}
\newcommand{\plus}{\oplus} \newcommand{\Plus}{\bigoplus}

\newcommand{\Q}{\ens Q}
\newcommand{\qq}{\forall\:}
\newcommand{\qu}{\call{Q}}

\newcommand{\res}[1]{{\left | {}_{#1} \right.}}

\newcommand{\rang}{\mathrm{rang}}

\newcommand{\rg}{\mathrm{rg}}

\newcommand{\saut}{\vskip 2em}

\newcommand{\si}{\mbox{ si }} 
\newcommand{\sinon}{\mbox{ sinon}}

  \newcommand{\sta}{\stackrel}
\newcommand{\sub}{\subseteq}

 \newcommand{\tens}{\otimes} 

\newcommand{\tenso}[1]{\raisebox{-1.5ex}{$\ma \stackrel{\displaystyle
\tens}{\scriptstyle #1}$}} 

\newcommand{\tilda}{\widetilde} 

\newcommand{\tq}{\: : \:}

\newcommand{\uni}{\bigcup}

\newcommand{\vide}{\emptyset}

\newcommand{\Z}{\ens Z} 
\newcommand{\ZZ}{\mbox{\Large $\mathbbm Z$}}
\newcommand{\zzz}{{\mathbf{z}}}

\makeindex

\begin{document}
\begin{center}
{
\bf \Large 

COHOMOLOGIE DES FIBRÉS EN DROITES SUR LES VARIÉTÉS MAGNIFIQUES DE RANG MINIMAL}

\end{center}

\begin{center}
{\large Alexis Tchoudjem}\\

Institut Camille Jordan\\
Universit\'e Claude Bernard Lyon I\\
Boulevard du Onze Novembre 1918\\
69622 Villeurbanne\\
FRANCE\\
Alexis.Tchoudjem@math.univ-lyon1.fr
\end{center}
\vskip 1cm 

{\bf R\'esum\'e} : Le  th\'eorème de Borel-Weil-Bott d\'ecrit la cohomologie des fibr\'es en droites sur les vari\'et\'es de drapeaux. On g\'en\'eralise ici ce th\'eorème \`a une plus grande classe de vari\'et\'es projectives : les vari\'et\'es magnifiques de rang minimal.

{\bf Abstract :}{\it (Cohomology of line bundles over wonderful varieties of minimal rank)} 
The Borel-Weil-Bott theorem describes the cohomology of line bundles over flag varieties. Here, one generalizes this theorem to a wider class of projective varieties : the wonderful varieties of minimal rank.
\tableofcontents


\mainmatter

\newpage
\section*{Introduction}

On se place sur un corps $\kk$ alg\'ebriquement clos de caract\'eristique nulle et on considère un groupe lin\'eaire $G$ semi-simple et connexe sur $\kk$.

On se donne une vari\'et\'e projective $X$ munie d'une action de $G$ et un fibr\'e en droites\infra{
en fait, dans l'article on raisonnera  plutôt avec des faisceaux inversibles
} $\pi : L \to X$. On suppose que $\pi$ est $G-$lin\'earis\'e \ie : $G$ agit sur $L$, $\pi$ est $G-$\'equivariant et l'action de $G$ est lin\'eaire dans les fibres de $\pi$.

La suite des groupes de cohomologie $H^d(X,L)$ ($d \in \Z_{\ge 0}$) forme alors une suite de repr\'esentations de dimension finie du groupe $G$. 

Quelles sont ces repr\'esentations du groupe $G$  ?

En degr\'e $d=0$, \cad pour l'espace des sections globales, on trouve dans \cite{B89} une description complète  pour le cas o\`u la vari\'et\'e $X$ est {\it sph\'erique}\index{sph\'erique (vari\'et\'e)} (\ie normale et avec une orbite ouverte d'un sous-groupe de Borel de $G$).  Mais, en degr\'e $d$ quelconque, il n'y a pas de r\'eponse dans un cadre aussi vaste.

N\'eanmoins, lorsque $X$ est homogène, autrement dit une vari\'et\'e de drapeaux, il y a le th\'eorème de Borel-Weil-Bott qui d\'ecrit très simplement les groupes de cohomologie $H^d(X,L)$ : ils sont tous nuls sauf au plus en un degr\'e, o\`u l'on obtient une repr\'esentation irr\'eductible de $G$.

On a aussi une description explicite dans le cas des compactifications de groupes (\ie des compactifications de l'espace homogène $K \croi K / K$ pour un groupe semi-simple $K$), \cf \cite{K} et \cite{toi}.

Le but de cet article est de g\'en\'eraliser le th\'eorème de Borel-Weil-Bott \`a la classe des {\it vari\'et\'es magnifiques de rang minimal}. Les vari\'et\'es magnifiques ont \'et\'e introduites dans \cite{DCP} et \cite{L}. D'après \cite{L}, les vari\'et\'es magnifiques sont toutes sph\'eriques. Les vari\'et\'es magnifiques de rang minimal sont particulièrement \'etudi\'ees dans \cite{Nico}. Cette classe de vari\'et\'es, dont nous rappellerons la d\'efinition, comprend notamment les vari\'et\'es de drapeaux et les compactifications magnifiques au sens de \cite{DCP} des espaces homogènes $K \croi K / K$ ($K$ est un groupe adjoint), $PGL_{2n}/PSp_{2n}$, $E_6/F_4$. 

Afin d'obtenir la description des groupes de cohomologie des fibr\'es en droites sur une vari\'et\'e $X$ magnifique et de rang minimal, on utilise une d\'ecomposition de $X$ en cellules de Bialynicki-Birula. On fait ensuite intervenir un complexe de Grothendieck-Cousin qui met en jeu des groupes de cohomologie \`a support. Si on note $\goth g$ l'algèbre de Lie de $G$, ces groupes de cohomologie \`a support sont naturellement des $\goth g-$modules. L'\'etude de ces $\goth g-$modules et de la d\'ecomposition cellulaire de la vari\'et\'e $X$ se trouve simplifi\'ee  quand on se place dans le cadre des vari\'et\'es magnifiques de rang minimal. Et cela suffit pour arriver au r\'esultat. 

Cette m\'ethode est celle utilis\'ee dans \cite{K} et \cite{toi} pour obtenir, pour chaque groupe adjoint $K$, la cohomologie des fibr\'es en droites sur la compactification de l'espace homogène $K \croi K / K$. On ajoute ici l'\'etude de certaines courbes irr\'eductibles sur la vari\'et\'e $X$ (\cf les sections \ref{sec:chaine} et \ref{sec:etu}) et cela permet de traiter ensemble les cas de toutes les vari\'et\'es magnifiques de rang minimal.

Avant d'\'enoncer le th\'eorème principal (le th\'eorème \ref{thm:principal}) de cet article on va introduire quelques notations et rappeler quelques d\'efinitions :

\section{Notations concernant le groupe}
Soit $G$ un groupe alg\'ebrique lin\'eaire semi-simple connexe et simplement connexe sur $\kk$, d'algèbre de Lie $\goth g$. On choisit $B$ un de ses sous-groupes de Borel et $T$ un tore maximal de $B$ ; on appelle $B^-$ le sous-groupe de Borel oppos\'e \`a $B$, relativement \`a $T$ (\ie tel que $B^- \cap B = T$). Soient $\Phi$ et $W$ le syst\`eme de racines et le groupe de Weyl de $(G,T)$. Pour toute racine $\alpha \in \Phi$, on note $\alpha^\ch : \kk^* \to T$ la coracine correspondante. On notera $\Phi^+$  l'ensemble des racines positives relativement \`a $B$, $\rho$ la demi-somme des racines positives, et, si $w \in W$, on pose $w * \lambda := w(\lambda + \rho) - \rho$ pour tout caract\`ere $\lambda$ de $T$. Soient $\Delta$ la base de $\Phi$ d\'efinie par $B$ et  $l$ la fonction longueur correspondante sur $W$. 

Soit ${\mathcal X}$ le r\'eseau des caractères de $T$. \'Etant donn\'es un caractère $\lambda$ de $T$ et un sous-groupe \`a un paramètre $\nu : \kk^* \to T$ , on notera $\cg \lambda , \nu \cd$ l'unique entier tel que :

$$\qq s \in \kk^*, \lambda(\nu(s)) = s^{\cg \lambda , \nu \cd} \p$$

On dira qu'un caractère $\lambda$ de $T$ est {\it dominant}\index{dominant} si pour toute racine positive $\alpha \in \Phi^+$, $\cg \lambda,\alpha^\ch \cd \ge 0$ et qu'il est {\it r\'egulier}\index{r\'egulier} si pour toute racine positive $\alpha \in \Phi^+$, $\cg\lambda,\alpha^\ch \cd \not= 0$.

Soit $(\omega_\delta)_{\delta \in \Delta}$ la base des poids fondamentaux ; elle est form\'ee des caractères de $T$ qui v\'erifient  :
$$\qq \delta , \epsilon \in \Delta, \cg \omega_\delta , \epsilon^\ch \cd = \left\{ \begin{array}{cl}
1 & \mbox{ si } \delta = \epsilon \\
 0 & \sinon .
\end{array}
\right.$$ 

Lorsque $\lambda$ est un caractère de $T$ tel que $\lambda + \rho$ est r\'egulier, il existe un unique  $w_\lambda \in W$\index{$w_\lambda$} tel que $w_\lambda * \lambda$ est un caractère dominant. Dans ce cas, on note $\lambda^+ : = w_\lambda*\lambda$\index{$\lambda^+$} ce poids dominant et on pose : $l(\lambda) := l(w_\lambda)$\index{$l(\lambda)$}. 

Pour tout caractère dominant $\lambda$ de $T$, on note $L(\lambda)$\index{$L(\lambda)$, $G-$module simple de plus haut poids $\lambda$} le $G-$module irr\'eductible de plus haut poids $\lambda$.

 Enfin, on notera $( \cdot,\cdot)$ un produit scalaire $W-$invariant sur ${\mathcal X} \tenso{\Z}{\Q}$.

\section{Vari\'et\'es magnifiques}

Suivant \cite{L}, une vari\'et\'e alg\'ebrique $X$ munie d'une action du groupe $G$ est appel\'ee {\it magnifique}\index{magnifique} si :
\vskip .2cm
1) $X$ est lisse et compl\`ete ;

2) $G$ poss\`ede une orbite dense, $X^0_G$, dans $X$ dont le bord, $X \moins X^0_G$ est une r\'eunion de diviseurs premiers $D_i$, $i \in \{1 , ..., r\}$, lisses, \`a croisements normaux et d'intersection non vide ;

3) pour tous points $x,x' \in X$, si $\{i \tq x \in D_i\} = \{i \tq x' \in D_i\}$, alors $G.x = G.x'$.

\begin{center}
*
\end{center}

En fait, la $G-$ vari\'et\'e $X$ est entièrement d\'etermin\'ee par son orbite ouverte $X^0_G$. On dit aussi que $X$ est la {\it compactification magnifique} de l'espace homogène $X^0_G$.

Dans le point 2) de la d\'efinition les $D_i$ sont les {\it diviseurs limitrophes}\index{diviseur limitrophe} de $X$ et 
l'entier $r$ est le {\it rang}\index{rang} de la vari\'et\'e magnifique $X$. La vari\'et\'e magnifique $X$ a exactement $2^r$ orbites de $G$ dont une seule est ferm\'ee. Nous noterons $F$\index{$F$, unique $G-$orbite ferm\'ee de $X$} cette orbite. 

Par exemple, les vari\'et\'es magnifiques de rang $0$ sont les vari\'et\'es de drapeaux $G/P$ où $P$ est un sous-groupe parabolique de $G$.
Une autre famille de vari\'et\'es magnifiques est form\'ee par les compactifications magnifiques de groupes adjoints construites dans \cite[\S 3.4]{DCP}. Dans ce cas le rang $r$ est le rang du groupe.

Les vari\'et\'es de drapeaux et les compactifications de groupes font partie d'une classe commune de vari\'et\'es magnifiques : les {\it vari\'et\'es magnifiques de rang minimal}. Nous rappelons leur d\'efinition ci-dessous.

\subsection{Vari\'et\'es magnifiques de rang minimal}
 
Soit $X$ une $G-$vari\'et\'e magnifique.

L'orbite ouverte de $G$ dans $X$ est isomorphe \`a l'espace homogène $G/H$ o\`u $H$ est un sous-groupe ferm\'e de $G$. Le rang de $X$ v\'erifie toujours :

$$r \ge \rang(G) - \rang (H) $$

(o\`u le rang d'un groupe est la dimension de ses tores maximaux)

\begin{defin}
Lorsque $r = \rang(G) - \rang(H)$, on dit que $X$ est une {\it vari\'et\'e magnifique de rang minimal}\index{rang minimal}.
\end{defin}

N. Ressayre a classifi\'e toutes les vari\'et\'es magnifiques de rang minimal : elles s'obtiennent toutes par produit, par recouvrement fini ou par induction parabolique (\cf \cite[\S 3.4]{L01} ou la d\'efinition \ref{defi:indupara} p. \pageref{defi:indupara}) \`a partir des vari\'et\'es de drapeaux et des compactifications magnifiques des espaces homogènes suivants (\cf \cite{Nico}) :

$$K \croi K / K \; \mbox{ pour un groupe $K$ adjoint }\;;$$\label{liste}
$$ PGL_{2n} /PSp_{2n} \: , ; n \ge 2 \;;$$
$$E_6/F_4 \;;$$
$$PSO_{2n}/PSO_{2n-1}\;;$$
$$SO_7 /G_2 $$
(pour se ramener au cas simplement connexe, les compactifications de ces espaces homogènes sont munies de l'action du rev\^etement universel de $G$).

\subsection{Racines sph\'eriques}\label{sec:racsp}

Soit $X$ une $G-$vari\'et\'e magnifique. Comme pour le groupe $G$, on va fixer quelques notations concernant la vari\'et\'e $X$ : notamment un sous-groupe parabolique $Q_X$ de $G$ et un ensemble fini $\Sigma_X$ de caract\`eres de $T$ : les {\it racines sph\'eriques}\index{racines sph\'eriques} de $X$. 

Pour d\'efinir $Q_X$ et $\Sigma_X$, on note $\zzz$ le {\it point-base}\index{point-base} de $X$ ; c'est l'unique point fixe du sous-groupe de Borel $B^-$. Le sous-groupe parabolique $Q_X$\index{$Q_X=Q$, sous-groupe parabolique associ\'e \`a $X$} est le stabilisateur de $\zzz$ dans $G$.

Le point $\zzz$ est dans l'orbite ferm\'ee $F$ de $X$ et si on note $T_\zzz X$ et $T_\zzz F$ les espaces tangents en $\zzz$ de $X $ et de $F$, $\Sigma_X$\index{$\Sigma_X$, ensemble des racines sph\'eriques de $X$} est l'ensemble des poids du tore $T \sub Q_X$ dans le quotient $T_\zzz X / T_\zzz F$.

Dor\'enavant, on notera $Q=Q_X$ le sous-groupe parabolique associ\'e \`a $X$. 

{\bf Remarque :} en fait, l'ensemble des racines sph\'eriques, $\Sigma_X$, est une base d'un système de racines (\cf \cite[\S 3.4 et 3.5]{B90a})

\subsection{Groupe de Picard}

Soit $X$ une $G-$vari\'et\'e magnifique. 

Rappelons ici comment est associ\'e un caractère du tore \`a chaque faisceau inversible $\call{L}$ sur $X$.

Un tel faisceau admet une unique $G-$lin\'earisation sur $X$ car le groupe $G$ est simplement connexe (\cf \cite[lem. 2.2 et pro. 2.3]{KKV}).

Le faisceau $\call{L}$ \'etant $G-$lin\'earis\'e, le sous-groupe parabolique $Q$ opère sur la fibre $\call{L}\res{\zzz}$ (en le point-base $\zzz$) via un caractère $\lambda$ : on appelle ce caractère le {\it poids du faisceau} $\call{L}$ \index{poids d'un faisceau inversible}.

 On notera $$\pic (X) \sub {\mathcal X}$$
\index{$\pic(X)$} le sous-r\'eseau des caractères de $T$ qui sont le poids d'un faisceau inversible sur $X$.

\begin{center}
*
\end{center}

{\bf Remarques :}

\primo) Notons $\Pic (X)$ et $\Pic(F)$ les groupes de Picard de la vari\'et\'e $X$ et de son orbite ferm\'ee $F$. Si on note ${\mathcal X}_Q \sub {\mathcal X}$ le r\'eseau des caractères de $Q$, comme $F$ est isomorphe \`a la vari\'et\'e de drapeaux $G/Q$, on a l'isomorphisme :
$$\Pic (F) \iso {\mathcal X}_Q$$

et $\pic(X)$ est l'image de $\Pic (X)$ par la restriction :
$$\Pic(X) \to \Pic(F) \iso {\mathcal X}_Q \sub {\mathcal X} \p$$

De plus, d'après \cite[pro. 8.1]{DCP}, la restriction $\Pic(X) \to \Pic(F)$ et donc $\Pic(X) \to \pic(X)$ sont injectives, autrement dit un faisceau inversible sur $X$ est uniquement d\'etermin\'e par son poids.

\secundo) Si $H$ est le groupe d'isotropie d'un point de la $G-$orbite ouverte de $X$ alors, d'apr\`es \cite[lemme 2.2]{BLV}, $\pic(X)$ est le sous-r\'eseau de ${\mathcal X}$ engendr\'e par les caract\`eres dominants $\lambda$ tels que le $G-$module simple $L(\lambda)$ ait au moins un ${H}-$vecteur propre.

\begin{center}
*
\end{center}

\section{Cohomologie des fibr\'es en droites}
Maintenant que sont fix\'ees les notations concernant le groupe $G$, la vari\'et\'e magnifique $X$ et les faisceaux inversibles sur $X$, on va \'enoncer le th\'eorème principal de ce texte. Il d\'ecrit, pour toutes les vari\'et\'es magnifiques de rang minimal,  les groupes de cohomologie de tous les faisceaux inversibles.

\begin{thm}\label{thm:principal}

Soit $X$ une vari\'et\'e magnifique de rang minimal. En tout degr\'e $d \ge 0$ et pour tout faisceau inversible $\call{L}_\lambda$ sur $X$ de poids $\lambda \in \pic(X)$, on a un isomorphisme de $G-$modules :{\LARGE
$$H^d(X,\call{L}_\lambda) \iso \Plus_{J \sub \Sigma_X} \Plus_{{\mu \in (\lambda + R_J) \cap \Omega_J \atop
\mu + \rho \mbox{ \scriptsize r\'{e}gulier }} \atop  l(\mu) +  |J|= d } L(\mu^+) $$
}
o\`u $\Sigma_X$ est l'ensemble (fini) des racines sph\'eriques de $X$ et pour toute partie $J$ de $\Sigma_X$, 
$$R_J :=  \sum_{\gamma \in J} \ZZ_{>0} \gamma + \sum_{\gamma \in \Sigma_X \moins J} \ZZ_{\le 0} \gamma $$
$$ \et \Omega_J := \Bigm\{ \mu \in \pic(X) \tq  \{ \gamma \in \Sigma_X \tq  (\mu + \rho ,\gamma ) < 0 \} = J \Bigm\} \p$$

\end{thm} 
\begin{center}
*
\end{center}

{\bf Remarques :} 

Ce th\'eorème g\'en\'eralise un des r\'esultats de \cite{K} et \cite{toi} qui concerne les compactifications de groupes adjoints. On retrouve aussi quelques r\'esultats d\'ej\`a connus (dans le cadre plus g\'eneral des vari\'et\'es sph\'eriques) :

\begin{liste}
\item en degr\'e $0$ : 
$$H^0(X,\call{L}_\lambda) = \Plus_{\mu} L(\mu) $$
o\`u $\mu$ d\'ecrit les poids dominants de l'ensemble $\ma \lambda + \sum_{\gamma \in  \Sigma_X} \Z_{\le 0} \gamma\;$ (\cf \cite[pro. 2.4]{B89}) ; 

\item lorsque $\lambda$ est dominant tous les groupes de cohomologie sup\'erieure (\ie en degr\'e $d >0$) sont nuls (\cf \cite[\S 2.1, cor. 1]{B90b}) ;

\item les multiplicit\'es des $G-$modules irr\'eductibles qui apparaissent dans les groupes de cohomologie $H^d(X,\call{L}_\lambda)$ peuvent \^etre $> 1$ (\cf \cite[\S 2.2, rem. 6]{toi}) ; 

\item les multiplicit\'es  sont born\'ees par $|W|$, ind\'ependamment du faisceau $\call{L}_\lambda$ (en effet, si $\mu_0$ est un caractère dominant, les caractères $\mu$ tels que $\mu^+ = \mu_0$ sont tous dans l'orbite $W * \mu_0$) (\cf aussi \cite[th. 3.4]{B94}).
\end{liste}

\begin{center}
*
\end{center}
({\it Dualit\'e de Serre})
 
Notons  $\ma 2\rho_X := \!\!\!\!\sum_{\alpha >0 \atop \alpha \; \mathrm{non \; racine\; de\;}Q_X} \!\!\alpha$ et posons : $J^* := \Sigma_X \moins J$ et $\mu^* : = -\mu - 2 \rho_X$, pour toute partie $J$ de $\Sigma_X$ et tout $\mu \in \pic(X)$. On v\'erifie la dualit\'e de Serre gr\^ace \`a l'involution  :
$$(J , \mu) \donne (J^*, \mu^*) $$
de l'ensemble ${\mathcal P}(\Sigma_X) \croi \pic(X)$ \index{${\mathcal P}(\Sigma_X)$, l'ensemble des parties de $\Sigma_X$}.

\subsection{Cas particuliers}

\subsubsection{Vari\'et\'es de drapeaux}

Dans le cas des vari\'et\'es de drapeaux $X = G/Q$, o\`u $Q$ est un sous-groupe parabolique de $G$, il n'y a pas de racines sph\'eriques : $\Sigma_X= \vide$ et donc :
$$H^d(X, \call{L}_\lambda) = \left\{\begin{array}{cl}
L(\lambda^+) & \si l(\lambda) =d \;;\\
0 & \sinon .
\end{array}
\right.$$
On retrouve donc le th\'eorème de Borel-Weil-Bott.

\subsubsection{Vari\'et\'es sym\'etriques complètes de rang minimal}\index{sym\'etrique complète (vari\'et\'e)}

Dans ce paragraphe, on suppose que le groupe $G$ est adjoint (\ie de centre trivial).

\'Etant donn\'ee une involution $\theta : G \to G$, on note $G^\theta$ le sous-groupe :
$$\left\{ g \in G \tq \theta(g) = g \right\} \p$$

Chaque espace sym\'etrique $G/G^\theta$ admet une unique compactification magnifique $X$ (\cf \cite{DCP}), appel\'ee une {\it vari\'et\'e sym\'etrique complète}.

Les vari\'et\'es sym\'etriques complètes sont des $\tilda{G}-$vari\'et\'es magnifiques, pour le rev\^etement universel $\tilda{G}$ de $G$ ; et leur rang est minimal uniquement pour les couples $(G ,G^\theta)$ suivants (qui apparaissent dans la liste p. \pageref{liste}) :
$$(K \croi K , \diag (K)) \;, \; (PGL_{2n},PSp_{2n})\;,\; (PSO_{2n}, PSO_{2n-1})\;,\; (E_6,F_4)\;$$
(o\`u $K$ est un groupe adjoint) \cf \cite{Nico}.

Notons respectivement :
$$\adh{K} \;,\; \adh{PGL_{2n}/PSp_{2n}} \;,\; \adh{PSO_{2n}/PSO_{2n-1}} \;,\; \adh{E_6/F_4}$$
les vari\'et\'es magnifiques correspondantes.

\begin{center}
*
\end{center}

Le tableau suivant concerne l'annulation des groupes de cohomologie $H^d(X,\call{L}_\lambda)$ 
pour tous les faisceaux inversibles $\call{L}_\lambda$ sur une vari\'et\'e sym\'etrique complète $X$ de rang minimal et de dimension $N$\infra{
Les vari\'et\'es $\adh{K},\adh{PGL_{2n}/PSp_{2n}},\adh{PSO_{2n}/PSO_{2n-1}},\adh{E_6/F_4}$ sont respectivement de dimension : $N=\dim K, 2n^2-n-1,2n-1, 26$.} :
$$
\begin{array}{|l|l|c|}
\hline
X=\adh{K} &  d \ou N-d = 1,2, 4 &\\
\cline{1-2}
X = \adh{PGL_{2n}/PSp_{2n}} & d \ou N-d = 1,2,3,4,6,7, 8 & \multirow{2}{7em}{$H^d(X,\call{L}_\lambda) = 0$}\\
\cline{1-2}
X=\adh{PSO_{2n}/PSO_{2n-1}} & d \not= 0,2n-1 &\\
\cline{1-2}
X=\adh{E_6/F_4} &  d \not= 0 , 9, 17,26 &\\
\hline
\end{array}
$$

{\bf Remarque : }La $3$ème ligne du tableau est bien connue puisque la vari\'et\'e $\adh{PSO_{2n}/PSO_{2n-1}}$ est l'espace projectif $\PP^{2n-1}$. 

\begin{center}
*
\end{center}

Pour obtenir ces annulations, on montre que certains entiers $d$ ne sont jamais de la forme $l(\mu) + |J|$ pour une partie $J$ de $\Sigma_X$ et un  poids $\mu \in \Omega_J$ :

Soit $T_1$ un tore {\it $\theta-$anisotrope maximal}\index{tore anisotrope} \ie : 
$$\qq x \in T_1, \theta(x) = x\inv$$
et $T_1$ est maximal pour cette propri\'et\'e.

Suivant \cite[\S 1.1]{DCP}, on choisit un tore maximal $T$ de $G$ qui contient $T_1$ ; forc\'ement, $T$ est invariant par l'involution $\theta$ et si on note $T_0$ la composante neutre de $T^\theta := \{x \in T \tq \theta(x) = x\}$, on a :
\begin{equation}\label{eq:ttt}
T= T_0.T_1 \p
\end{equation}

On note encore $\theta$ l'involution induite sur le $\Q-$espace vectoriel ${\mathcal X}\tenso{\Z} \Q $. Lorsque $X$ est une vari\'et\'e sym\'etrique complète de rang minimal, on a :
$$\pic(X)\tenso{\Z}\Q = \left\{ \lambda \in {\mathcal X}\tenso{\Z} \Q \tq \theta(\lambda) = -\lambda\right\} $$
(\cf \cite[\S 8.1, rem.]{DCP}).

On peut aussi choisir l'ensemble des racines $\Phi^+$ tel que :
$$\qq \alpha \in \Phi^+, \theta(\alpha) = \alpha \mbox{ ou } \theta(\alpha) \in - \Phi^+ $$(\cf \cite[lem. 1.2]{DCP})

On pose alors $\Phi_0 := \{\alpha \in \Phi \tq \theta(\alpha) = \alpha\}$, $\Phi_1 := \Phi \moins \Phi_0$  et $\Phi_1^+ := \Phi_1 \cap \Phi^+$. On note $\Delta$ la base correspondant \`a $\Phi^+$. On a alors :

$$\Sigma_X=\left\{ \alpha - \theta(\alpha) \tq \alpha \in \Delta \cap \Phi_1^+ \right\}\p$$

 On pose ensuite $\tilda{\alpha} := \alpha - \theta(\alpha)$ pour toute racine $\alpha$ et :
$$\widetilde{\Phi} := \left\{ \alpha - \theta(\alpha) \tq \alpha \in  \Phi_1 \right\}
 \p$$

En fait, $\widetilde{\Phi}$ est un système de racines (non forc\'ement r\'eduit) dont $\Sigma_X$ est une base. L'ensemble des racines positives correspondant est :
$$\widetilde{\Phi}^+ := \left\{ \alpha - \theta(\alpha) \tq \alpha \in  \Phi_1^+\right\}
 \p$$ 

Soit maintenant $\mu \in \pic(X)$ ; du fait que : $\theta(\mu) = -\mu$, on d\'eduit que :
$$\{\beta \in \Phi^+ \tq (\mu + \rho , \beta) < 0 \} = \{\beta \in \Phi_1^+ \tq (\mu + \rho, \beta) <0\} \p$$

Or, si $\beta \in \Phi_1^+$, on a les \'equivalences :
$$(\mu + \rho , \beta) <0 \equi (\mu + \rho , -\theta(\beta) ) < 0 \equi (\mu + \rho , \tilda{\beta}) < 0 \p$$

Soit $\ma \widetilde{l}(\mu) := \left| \{ \tilda{\beta} \in \widetilde{\Phi}^+ \tq (\mu + \rho, \tilda{\beta)} < 0 \} \right| \p$

Lorsque $\tilda{\beta}$ d\'ecrit l'ensemble $\tilda{\Phi}^+$,  les ensembles :
$$\{ \alpha \in \Phi_1^+ \tq \tilda{\alpha} = \tilda{\beta} \}$$
restent de m\^eme cardinal. En effet, soient $\beta_0,\beta_1 \in \Phi^+$ ; il r\'esulte de  (\ref{eq:ttt}) que :
$$\tilda{\beta_0} = \tilda{\beta_1} \equi \beta_0 \res{T_1} = \beta_1\res{T_1} \p$$
Mais d'après \cite[pro. 4.7]{Ri}, il existe $n \in N_K(T_1)$ et $w \in W$ tels que :
$$n. (\beta_0 \res{T_1}) = \beta_1\res{T_1} \et \qq x \in T_1, w\inv x w = n\inv x n\p$$
Il s'ensuit que :
$$w.\{\alpha \in \Phi_1^+ \tq \tilda{\alpha} = \tilda{\beta_0} \} = \{\alpha \in \Phi_1^+ \tq \tilda{\alpha} = \tilda{\beta_1} \} \p$$

On v\'erifie dans chaque cas : 
$$X = \adh{K},\adh{PGL_{2n}/PSp_{2n}},\adh{PSO_{2n}/PSO_{2n-1}} \ou \adh{E_6/F_4}$$ que le cardinal commun :
$$|\{ \alpha \in \Phi_ 1^+ \tq \tilda{\alpha} = \tilda{\beta} \}|$$
($\tilda{\beta} \in \tilda{\Phi}^+$) vaut respectivement : $2, 4, (2n-2) \ou 8$.

Il en r\'esulte pour toute partie $J$ de $\Sigma_X$ que :
$$l(\mu) + |J| = \left\{ \begin{array}{cc}
2\tilda{l}(\mu) + |J| & \si X = \adh{K} ,\\
4\tilda{l}(\mu) + |J| &\si X = \adh{PGL_{2n}/PSp_{2n}},\\
(2n-2) \tilda{l}(\mu) + |J| &\si X =\adh{PSO_{2n}/PSO_{2n-1}},\\
8\tilda{l}(\mu) + |J| &\si X = \adh{E_6/F_4} \p
\end{array}\right.$$

Enfin on utilise que pour une partie $J$ de $\Sigma_X$ et pour un $\mu \in \Omega_J$, on a les in\'egalit\'es :
$$\tilda{l}(\mu) \ge |J| \et | \tilda{\Phi}^+ | - \tilda{l}(\mu) \ge |\Sigma_X \moins J|$$
et les \'equivalences :
$$\tilda{l}(\mu) = 0 \equi |J| = 0 \et  \tilda{l}(\mu) = |\tilda{\Phi}^+| \equi |J| = |\Sigma_X| \p$$
\hfill{\underline{\bf Q.e.d.}}  

\begin{center}
***
\end{center}

\subsection{Quelques figures en rang petit }\label{sec:dessins}

On traite ici le cas des vari\'et\'es magnifiques $X$, de rang $\le 2$,  suivantes :

\begin{liste}
\item (de rang $1$) l'espace projectif $\PP^{2n-1}$ vu comme compactification magnifique de l'espace homogène $PSO_{2n} / PSO_{2n-1}$ et $\PP^7$ vu comme compactification magnifique de l'espace homogène $SO_7/G_2$.

\item (de rang $2$) les compactifications magnifiques des vari\'et\'es sym\'etriques : $PGL_3 \croi PGL_3 / PGL_3$, $PGL_6/PSp_6$ et $E_6/F_4$.
\end{liste}

Sur les figures qui suivent, sont repr\'esent\'es les poids des ensembles $\Omega_J$ et $R_J$ de l'\'enonc\'e du th\'eorème principal avec les notations suivantes (\cite{Was}) :

\begin{liste}
\item pour les cas de rang $1$, on a not\'e $\gamma$ la racine sph\'erique de $X$ et $\widetilde{\omega}$ le g\'en\'erateur de $\pic(X)$ tel que $(\tilda{\omega} , \gamma) >0$. On a pos\'e  $\lambda_0 := -\left( \frac{(\rho,\gamma)}{(\tilda{\omega} , \gamma)} +1 \right)  \tilda{\omega} \in \pic(X)$ de sorte que si $n \in \Z$ on ait :
$$\lambda_0 +n \tilda{\omega} \in \Omega_{\Sigma_X} \equi n  \le 0 \;\;;$$

\item pour les cas de rang $2$, on a not\'e $\gamma_1,\gamma_2$ les racines sph\'eriques de $X$ et  $\widetilde{\omega_1}, \widetilde{\omega_2}$ la base de $\pic(X)$ telle que $(\tilda{\omega_i} , \gamma_i) >0 \et (\tilda{\omega_i} , \gamma_j) =0 \si i\not=j $ sinon. Enfin, on a pos\'e   $\lambda_0 := -\left( \frac{(\rho,\gamma_1)}{(\tilda{\omega_1} , \gamma_1)} +1 \right)  \tilda{\omega_1} - \left( \frac{(\rho,\gamma_2)}{(\tilda{\omega_2}, \gamma_2)} +1 \right) \tilda{\omega_2}  \in \pic(X)$ de sorte que si $n_1,n_2 \in \Z$, on ait :
$$\lambda_0 + n_1 \tilda{\omega_1} + n_2 \tilda{\omega_2} \in \Omega_J \equi  J = \{\gamma_i \tq n_i \le 0\} \p$$ 
\end{liste}

Pour les compactifications de $$PSO_{2n}/PSO_{2n-1} ,\;SO_7/G_2 ,\; PGL_3 \croi PGL_3 / PGL_3 ,\; PGL_6/PSp_6 \et E_6/F_4,\;$$  $\lambda_0$ 
d\'esigne respectivement les poids :
$$-n\widetilde{\omega} \;,\; -4 \widetilde{\omega}\;,\;
-2\widetilde{\omega_1} -2 \widetilde{\omega_2}\;,\; -3\widetilde{\omega_1} - 3 \widetilde{\omega_2} \;,\; -5\widetilde{\omega_1} -5 \widetilde{\omega_2} \p$$

\newpage

\begin{figure}[h]
    \input{b3g2P.pstex_t}
    \caption{Les ensembles de poids $\Omega_J$ pour les compactifications magnifiques de $PSO_{2n}/PSO_{2n-1}$ et $SO_7/G_2$}
  \end{figure}

\begin{figure}[h]
    \input{b3g2R.pstex_t}
    \caption{Les ensembles de poids $R_J$ pour les compactifications magnifiques de $PSO_{2n}/PSO_{2n-1}$ et $SO_7/G_2$}
  \end{figure}

\newpage

\begin{figure}[h]
   \input{pgl3pgl3P.pstex_t}
\caption{Les ensembles de poids $\Omega_J$ pour les compactifications magnifiques de $PGL_3 \croi PGL_3 / PGL_3$, $PGL_6/PSp_6$ et $E_6/F_4$}
\end{figure}
\newpage
\begin{figure}[h]
    \input{pgl3pgl3R.pstex_t}
    \caption{Les ensembles de poids $R_J$ pour les compactifications magnifiques de $PGL_3 \croi PGL_3 / PGL_3$, $PGL_6/PSp_6$ et $E_6/F_4$}
  \end{figure}

\newpage

On a not\'e sur chaque figure par des $\bullet$ les poids de $\Omega_{\Sigma_X}$ et $R_{\Sigma_X}$, par des $\circ$ les poids de $\Omega_\emptyset$ et $R_\emptyset$, par des $+$ les poids de $\Omega_{\gamma_1}$ et $R_{\gamma_1}$ et par des $\cdot$ les poids de $\Omega_{\gamma_2}$ et $R_{\gamma_2}$.

\begin{center}
***
\end{center}

\begin{center} 
{\it Fixons maintenant et jusqu'\`a la fin une $G-$vari\'et\'e magnifique $X$ et un faisceau inversible $\call{L}_\lambda$ de poids $\lambda$.}
\end{center}

\begin{center}
*
\end{center}

Les grandes \'etapes de la d\'emonstration du th\'eorème principal, sont, d'abord, la d\'ecomposition de la vari\'et\'e $X$ en cellules de Bialynicki-Birula (\cf la section \ref{subsec:cel}) et des cellules en orbites du Borel, ensuite, le calcul de groupes de cohomologie \`a support dans les cellules et dans les orbites du Borel, et enfin, pour passer, de la cohomologie \`a support \`a la cohomologie usuelle, la d\'ecomposition d'un complexe : le complexe de Grothendieck-Cousin (\cf la partie \ref{sec:GC}).

Au cours de cette d\'emonstration, on utilisera notamment deux propri\'et\'es importantes des vari\'et\'es magnifiques de rang minimal :

\noindent --- {\bf propri\'et\'e 1 :} {\it La vari\'et\'e $X$ n'a qu'un nombre fini de points fixes du tore $T$.}

En effet, $X$ n'a qu'un nombre fini de $G-$orbites et chaque $G-$orbite contient un nombre fini de points fixes.

\noindent --- {\bf propri\'et\'e 2 :} {\it La vari\'et\'e $X$ n'a qu'un nombre fini de courbes $T-$invariantes} (\cf la partie \ref{sec:etu} et le lemme \ref{lem:clambda}).
  
Mais avant cela, nous allons \'etablir  un r\'esultat
sur les racines sph\'eriques des vari\'et\'es magnifiques de rang minimal.

\section{Racines et racines sph\'eriques}\label{sec:rac-sph}

Il peut arriver qu'une racine sph\'erique d'une vari\'et\'e magnifique soit un multiple d'une racine de $(G,T)$. Nous allons montrer, dans cette section, que cela est impossible lorsque la vari\'et\'e est de rang minimal. En effet, on a d'abord : 

\begin{pro}\label{pro:not}
Il existe un point $y$ de la $G-$orbite ouverte de $X$, de groupe d'isotropie $G_y =: H$ et $P$ un sous-groupe parabolique de $G$ tel que :
\begin{liste}
\item[1)] le groupe $N_G(H) / H$\index{
$N_G(H)$ le normalisateur de $H$ dans $G$
} est fini ;
\item[2)] $H \sub P$ et $R_u(H) = R_u(P)$ ;
\item[3)] $T \sub P$ ;
\item[4)] l'orbite $B.y$ est ouverte dans $X$ ;
\item[5)] la composante neutre $(T\cap H)^\circ$ est un tore maximal de $H$.
\end{liste}
\end{pro}

\begin{dem}
Soient $y$ un point de la $G-$orbite ouverte de $X$ et $H$ son groupe d'isotropie.

Comme $X$ est magnifique, d'après \cite[cor. 5.3]{BP}, le quotient $N_G(H)/H$ est fini. 

D'un autre c\^ot\'e, d'après \cite[cor. 3.9]{BT}, il existe un sous-groupe parabolique $P$ de $G$ tel que :
$$H \sub P \et R_u(H)\sub R_u(P) \p$$
  
Comme $X$ est aussi de rang minimal, d'après \cite{Nico}, on peut choisir $P$ pour que :
$$R_u(H) = R_u(P) \p$$

Ainsi, on a obtenu les propri\'et\'es 1) et 2). Mais a priori $P$ ne contient pas $T$. 

Or, il existe $g \in G$ tel que le sous-groupe parabolique $gPg\inv$ contienne $B^-$.

En remplaçant $y$ par $g.y$, $H$ par $gHg\inv$, $P$ par $gPg\inv$, 1),2),3) sont v\'erifi\'ees et $PB$ est ouvert dans $G$.

Puisque la vari\'et\'e homogène $G/H$ est sph\'erique de rang minimal, il existe $x \in G/B$ fix\'e par un tore maximal de $H$, appelons-le $T_H$, tel que :
$$H.x \mbox{ est ouvert dans }G/B \p$$

Puisque $PB/B$ est un ouvert de $G/B$ stable par $H$, $x$ est de la forme : $x = p\inv B/B$ pour un $p \in P$.

Cette fois, en remplaçant $y$ par $p.y$, $H$ par $pHp\inv$, $T_H$ par $pT_Hp\inv$ et en gardant $P$, 1), 2), 3), 4) sont v\'erifi\'ees avec de plus $T_H \sub P \cap B$.

Il existe alors $a \in P \cap B$ tel que $a T_H a\inv \sub T$.

Finalement, en remplaçant $y$ par $a.y$, $H$ par $aHa\inv$ et en gardant $P$, 1), 2) 3), 4), 5) sont v\'erifi\'ees. 

\end{dem}

On garde $H$ et $T_H = (T \cap H )^\circ$ comme dans l'\'enonc\'e de la proposition ci-dessus. En particulier, $T_H$ est non trivial (sauf si $G$ est trivial). On va d\'emontrer :

\begin{lem}\label{lem:racnonsph}
$i)$ Si $\alpha$ est une racine de $(G,T)$, alors $\alpha\res{T_H} \not= 0$; 

$ii)$ Si $\gamma$ est une racine sph\'erique de $X$, relativement \`a $B$, alors : $\gamma\res{T_H} = 0$.
\end{lem}

Il r\'esulte imm\'ediatement de ce lemme qu'aucune racine n'est combinaison lin\'eaire de racines sph\'eriques.

\begin{dem}[du lemme]

\begin{center}
{\LARGE \it i)}
\end{center}
{\bf 1er cas : si $H$ est r\'eductif :}

Dans ce cas, $T_H = C_{H^\circ}(T_H)$ le centralisateur du tore $T_H$ dans la composante neutre $H^\circ$ de $H$.

Donc, pour tout $x \in G/B$, d'après \cite[th. A]{Richardson}, le tore $T_H$ agit transitivement sur les composantes connexes de l'ensemble de points fixes :
$$(H^\circ x B/B)^{T_H} \p$$
L'ensemble $(H^\circ x B/B)^{T_H}$ est par cons\'equent fini.

Or, puisque $H$ est un sous-groupe de $G$ sph\'erique et de rang minimal, il existe $x \in (G/B)^{T_H}$ tel que $H^\circ x B/B$ est ouvert dans $G/B$. 

L'ensemble fini $(H^\circ x B/B)^{T_H}$ est alors un ouvert non vide de $(G/B)^{T_H}$. On en d\'eduit qu'au moins une composante connexe de $(G/B)^{T_H}$ est un point. 

Ce point est forc\'ement laiss\'e fixe par le centralisateur $C_G(T_H)$ qui est connexe. Le groupe $C_G(T_H)$ est en cons\'equence r\'esoluble. Mais cela n'est possible que si, pour toute racine $\alpha$ de $(G,T)$, $\alpha\res{T_H} \not= 0$.

{\bf 2ème cas : si $H$ est quelconque :} 

On choisit d'abord un sous-groupe parabolique $P$ de $G$ tel que :
$$T \sub P , H \sub P, R_u(H) = R_u(P) $$
(\cf la proposition \ref{pro:not}). On pose ensuite :
$$L := P/ R_u(P) \et K:= H/R_u(P) \p$$
L'espace homogène $L/K$ est encore sph\'erique et de rang minimal mais de plus, $L$ et $K$ sont r\'eductifs.

Le radical de $L$, $R(L)$, est son centre connexe donc :
$R(L) \sub N_G(H) /R_u(P)$. Remarquons que l'on peut consid\'erer $R(L)$ comme un sous-tore de $T$ vu que $T \cap R_u(P) = \{1\}$ ; on a donc aussi $(R(L) \cap H)^\circ \sub T_H / R_u(P) $.

Or, le quotient $N_G(H)/H$ est fini, donc $R(L)/(R(L) \cap H)^\circ$ est fini. 

\'Etant donn\'ee une racine $\alpha$ de $(G,T)$, deux possibilit\'es se pr\'esentent :
\begin{liste}
\item soit $\alpha \res{R(L)} \not= 0$, 

mais alors, $R(L)$ \'etant connexe, avec $N:=\left | R(L)/(R(L) \cap H)^\circ \right |$, on trouve :
$$\alpha\res{R(L)} \not=0 \impliq  N\alpha\res{R(L)} \not=0 $$
$$\impliq \alpha\res{(R(L) \cap H)^\circ} \not =0$$
$$\impliq \alpha \res{T_H} \not= 0 \;\;\;\; ;$$

\item soit $\alpha\res{R(L)} = 0$,

 mais alors $\alpha$ est en fait une racine de $(L,T)$ et donc $\alpha\res{T_H} \not=0$ gr\^ace au premier cas trait\'e. 
\end{liste} 
\vskip -1cm
\begin{center}
{ \LARGE \it ii)}
\end{center}
\vskip -.3cm
Une racine sph\'erique $\gamma$ de $X$, relativement \`a $B$, est en particulier un poids de $B$ op\'erant dans l'espace $\kk(X)$ des fonctions rationnelles sur $X$ (\cf \cite[\S 1.3]{L01}). Autrement dit, il existe une fonction $f$, non nulle, r\'egulière sur un ouvert $B-$ stable, $V$, de $X$ telle que :
$$\qq b \in B, \qq v \in V, f(b\inv . v) = \gamma(b) f(v) \p$$
En choisissant $v = y$ comme dans l'\'enonc\'e de la proposition \ref{pro:not}, on a $f(y) \not= 0$ et :
$$\qq t \in T_H, \; f(y) = f(t\inv.y) = \gamma(t) f(y)$$
$$\impliq \qq t \in T_H, \gamma(t) = 1 $$
$$\impliq \gamma\res{T_H} = 0\p$$\end{dem}

{\bf Remarque :} Dans le i), lorsque $H$ est semi-simple, pour toute racine $\alpha$ de $(G,T)$, $\alpha\res{T_H}$ est en fait une racine de $(H,T_H)$ (\cf \cite[lem. 4.2]{Nico}).

\section{D\'ecomposition cellulaire}

Du fait que la vari\'et\'e $X$ est projective et n'a qu'un nombre fini de points fixes du tore $T$,  on d\'eduit que $X$ s'\'ecrit comme une r\'eunion disjointe finie d'espaces affines centr\'es en les points fixes du tore $T$, les cellules de Bialynicki-Birula :

\subsection{Cellules}\label{subsec:cel}
Soit $x \in X^T$ un point fixe du tore $T$.

\'Etant donn\'e un sous-groupe \`a un param\`etre $\zeta : \kk^* \to T$, on notera :
 
$$X^+(x) \index{$X^+(x)$, cellule centr\'ee en $x$}:= \{y \in X \tq \lim_{a \to 0} \zeta(a) y = x \} $$
c'est la {\it cellule centr\'ee en $x$}.

D'apr\`es \cite{BiBi}, $X^+(x)$ est une sous-vari\'et\'e localement ferm\'ee de $X$ qui est isomorphe \`a un espace affine.

Plus pr\'ecis\'ement, l'espace tangent en $x$, $T_x X$ est un $T-$module et donc un $\kk^*-$module via $\zeta$ ; il en r\'esulte une d\'ecomposition en $\kk^*-$espaces propres :
$$T_x X = \Plus _{n \in \Z} (T_x X)_n$$avec $(T_x X)_n : = \left\{ v \in T_xX \tq \qq a \in \kk^* , \, \zeta(a) .v = a^n v \right\}$. On pose alors :
$$(T_xX)_+\index{$(T_xX)_+$, partie positive de l'espace tangent en $x$} := \Plus_{n >0} (T_xX)_n \et (T_x X)_- \index{$(T_x X)_-$, partie n\'egative de l'espace tangent en $x$}:= \Plus_{n < 0} (T_x X)_n \p$$

Avec ces notations, la sous-vari\'et\'e $X^+(x)$ est $T-$isomorphe \`a l'espace affine $(T_x X)_+$.

On utilisera surtout des sous-groupes \`a un paramètre qui ont le moins de points fixes possibles :

\begin{defin} 
On dira d'un sous-groupe \`a un param\`etre $\zeta$ qu'il est {\it  $X-$r\'egulier \index{X-r\'egulier }} si :
$$X^{\zeta(\kk^*)} = X^T $$
et qu'il est {\it dominant}\index{dominant} si pour toute racine positive $\alpha$, $\cg \alpha,\zeta \cd > 0 \p$ 
\end{defin} 

{\bf Remarques : } 

\primo) Dans le r\'eseau $\Hom(\kk^* , T)$, le compl\'ementaire de l'ensemble des sous-groupes \`a un paramètre $X-$r\'eguliers est une r\'eunion de sous-groupes stricts. En particulier, il existe toujours des sous-groupes \`a un paramètre $X-$r\'eguliers et dominants.

\secundo) si, par exemple, $X$ est la vari\'et\'e de drapeaux $ G/B^-$, alors $\zeta$ est un sous-groupe \`a un paramètre $X-$r\'egulier si et seulement si $\cg \alpha , \zeta \cd \not = 0$
pour toute racine $\alpha \in \Phi$.

\saut

Si $\zeta$ est un sous-groupe \`a un param\`etre dominant et $X-$r\'egulier, on obtient une d\'ecomposition cellulaire finie de $X$  :
$$X = \Dij_{x \in X^T} X^+(x)$$
dont les {\it cellules}\index{cellules} $X^+(x)$ sont $B-$invariantes (car si $\ma b \in B,$ alors :
$$ \lim_{a \to 0} \zeta(a)b\zeta(a)\inv \in T \;).$$ 

Remarquons que cette d\'ecomposition d\'epend en g\'en\'eral du choix du sous-groupe \`a un param\`etre dominant et $X-$r\'egulier $\zeta$. 
 
On va param\'etrer l'ensemble $X^T$ :

\subsection{L'ensemble des points fixes du tore}

Rappelons pour commencer une caract\'erisation remarquable des vari\'et\'es magnifiques de rang minimal :

\begin{pro}[{\cite[pro. 2.3]{Nico}}]\label{pro:pointsfixes}
Une $G-$vari\'et\'e magnifique $X$ est de rang minimal si et seulement si les points fixes du tore $T$ sont tous dans la $G-$orbite ferm\'ee \raisebox{-0.2ex}{$\Box$}
\end{pro}


\saut

On peut donc param\'etrer l'ensemble $X^T$ par une partie du groupe de Weyl $W$. En effet, soit $W_{Q}$\index{$W_{Q}$, groupe de Weyl de $Q$} le groupe de Weyl de $(Q,T)$ ; on pose :
$$W^{Q}\index{$W^{Q}$} := \{ w \in W \tq \qq v \in W_{Q} ,\; l(wv) = l(w) + l(v) \} \p$$

L'ensemble $W^{Q}$ est un syst\`eme de repr\'esentants du quotient $W/W_{Q}$ et param\`etre l'ensemble des points fixes de l'orbite ferm\'ee $F \iso G/Q$ :
$$X^T = \{w \zzz \tq w \in W^{Q}\}$$
o\`u $\zzz$ est l'unique point fixe de $Q$ dans $X$.

Pour tout $w \in W^Q$, on notera $X^+_w := X^+(w\zzz)$\index{$X^+_w$, cellule centr\'ee en $w\zzz$} la cellule centr\'ee en $w\zzz \in X^T$. Alors :

$$X = \Dij_{x \in X^T} X^+(x) = \Dij_{w \in W^Q}{X^+_w} \p$$

{\it D\'esormais, un sous-groupe \`a un paramètre $\zeta$, dominant et $X-$r\'egulier, est fix\'e et lorsqu'on consid\'erera  une cellule $C$, on sous-entendra que $C = X^+_w$ pour un certain $w \in W^Q$. En particulier les cellules consid\'er\'ees sont stables par $B$.}

\subsection{Codimension des cellules}

C'est le moment de rappeler (\cf \cite{L}) que l'ensemble $\Sigma_X$ des racines sph\'eriques de $X$ (\cf page \pageref{sec:racsp}) est en bijection avec l'ensemble des diviseurs limitrophes $D_i$ (les composantes irr\'eductibles de $X \moins X^0_G$).

Plus pr\'ecis\'ement, pour chaque $1 \le i \le r$, le poids $\gamma_i$ du faisceau inversible $G-$lin\'earis\'e ${\mathcal O}_X(D_i)$ est une des racines sph\'eriques ; on les obtient toutes ainsi et les $\gamma_i$ sont deux \`a deux distincts.

\begin{defin}\label{def:Dgamma}
On notera $D_\gamma$\index{$D_\gamma$, diviseur limitrophe associ\'e \`a la racine sph\'erique $\gamma$} le diviseur limitrophe de $X$ correspondant \`a la racine sph\'erique $\gamma$ (le tore maximal $T$ agit donc via le caractère $\gamma$ sur la fibre ${{\mathcal O}_X(D_\gamma)}_{\zzz}$).
\end{defin}

La proposition qui suit donne en particulier la codimension des cellules $X^+_w$ :
\begin{pro}\label{pro:coce}
Si $w \in W^Q$, alors :
\begin{liste}
\item[i)] pour toute racine sph\'erique $\gamma \in \Sigma_X$, $$X^+_w \sub D_\gamma \equi \cg w\gamma , \zeta \cd < 0 \;\;;
$$
\item[ii)] $$\codim_X X^+_w = l(w) + \left| \{\gamma \in \Sigma_X \tq \cg w\gamma , \zeta \cd <0 \}\right| \p$$
\end{liste}
\end{pro}

\begin{dem}
Posons $x:= w\zzz$ et fixons  une racine sph\'erique $\gamma_0 \in \Sigma_X$.

\begin{liste}
\item[i)] L'espace tangent en le point fixe $x \in X^T$ est un $T-$module qui se d\'ecompose :
$$T_xX = T_x D_{\gamma_0} \plus T_x X / T_x D_{\gamma_0}$$
$$= T_x D_{\gamma_0} \plus {{\mathcal O}_X(D_{\gamma_0})}_{x} \p$$

La $T-$droite propre ${{\mathcal O}_X(D_{\gamma_0})}_{x} = {{\mathcal O}_X(D_{\gamma_0})}_{w\zzz} $  a pour poids $w {\gamma_0}$ selon $T$ car ${{\mathcal O}_X(D_{\gamma_0})}_\zzz$ a pour poids $\gamma_0$ et le faisceau inversible ${{\mathcal O}_X(D_{\gamma_0})}$ est $G-$ lin\'earis\'e sur $X$.

Par cons\'equent :
$$(T_x X)_+ = \left\{ 
\begin{array}{ll}
(T_x D_{\gamma_0})_+ & \si \cg w {\gamma_0} , \zeta \cd < 0\\
(T_x D_{\gamma_0})_+ \plus {{\mathcal O}_X(D_{\gamma_0})}_{x} & \si \cg w {\gamma_0} , \zeta \cd > 0 \p
\end{array}
\right.
$$
Comme $X^+(x)$ est $T-$ isomorphe \`a l'espace affine $(T_x X)_+$, on en d\'eduit que :
$$X^+_w = X^+(x) \sub D_{\gamma_0} \equi \cg w {\gamma_0} , \zeta \cd < 0 \p$$

\item[ii)] La codimension de la cellule $X^+_w$ est donn\'ee par :
$$\codim_X X^+_w = \codim_X X^+(x) $$
$$= \dim T_x X - \dim (T_x X)_+ = \dim (T_x X)_- \p$$
Or :
$$(T_x X)_- = ( T_x F)_- \plus (T_x X / T_x F)_- \p$$

Mais d'une part : $(T_x F)_- = (T_{w\zzz} F)_- $ est un $T-$module de dimension $ l(w)$ et d'autre part, $T_x X /T_x F$ est un $T-$module dont les poids sont exactement les $w \gamma$, $\gamma$ d\'ecrivant l'ensemble $\Sigma_X$ des racines sph\'eriques de $X$. Donc :
$$\dim (T_x X)_- = l(w) +  \left| \{\gamma \in \Sigma_X \tq \cg w \gamma , \zeta \cd < 0 \} \right|\p$$
\end{liste}
\end{dem}

\saut
Pour calculer la cohomologie des fibr\'es en droites sur $X$, on va utiliser une stratification de $X$. Malheureusement, en g\'en\'eral, la d\'ecomposition cellulaire de $X$ n'est pas une stratification au sens o\`u l'adh\'erence d'une cellule n'est pas toujours une union d'autres cellules. C'est pourquoi on aura besoin d'une d\'ecomposition plus fine de la vari\'et\'e $X$.

\section{D\'ecomposition en orbites du sous-groupe de Borel}

\subsection{Filtration par une suite de ferm\'es embo\^it\'es}\label{par:dec}

On utilise la d\'ecomposition de $X$ en $B-$orbites. La vari\'et\'e $X$ ne possède qu'un nombre fini de $B-$orbites (\cf par exemple \cite{L}).
De plus, chaque $B-$orbite  est une sous-vari\'et\'e affine. 

On pose ensuite pour tout entier $i$ : 
$$Z_i := \uni_{ \call{B}  \; B-\mathrm{orbite} \atop \codim_X(\call{B}) \ge i } \call{B} \p$$
 Puisque le bord $\adh{\call{B}} \moins \call{B}$ d'une orbite $\call{B}$ est une r\'eunion (finie) d'orbites de codimension strictement sup\'erieure, les $Z_i$ sont des sous-vari\'et\'es ferm\'ees de $X$ telles que :
$$X = Z_0 \con Z_1 \con ...$$
et pour tout $i$, $Z_i \moins Z_{i+1}$ est une r\'eunion disjointe finie de sous-vari\'et\'es affines de m\^eme codimension : $i$.

\subsection{Lien entre les cellules et les orbites du sous-groupe de Borel}\label{subsec:cgo}

Comme les cellules de Bialynicki-Birula (pour un sous-groupe \`a un paramètre dominant et $X-$r\'egulier fix\'e) sont stables par $B$, chaque $B-$orbite est contenue dans une unique cellule.

De plus, on a une param\'etrisation des points fixes du tore et donc des cellules par $W^Q$. Il en r\'esulte aussi une param\'etrisation des $B-$orbites car d'après \cite[\S 2.1 p. 219 et pro. 2.3]{BL}, l'intersection d'une cellule et d'une $G-$orbite est soit vide soit une $B-$orbite ; on obtient ainsi, de façon unique, toutes les $B-$orbites de $X$. En r\'esum\'e, les $B-$orbites de $X$ sont les intersections non vides parmi les  
$$\call{O} \cap X^+_w \; , \; $$
o\`u $\call{O}$ est une $G-$orbite et $w \in W^Q$. 

\section{\'Etude de certains groupes de cohomologie \`a support}

Avant de d\'eterminer les groupes de cohomologie  $H^i(X,\call{L}_\lambda) $,
on s'int\'eresse d'abord, en guise d'approximation, aux groupes de cohomologie \`a support :
$$H^i_\call{B}(\call{L}_\lambda) \et H^i_C(\call{L}_\lambda)$$
dans une $B-$orbite $\call{B}$ et une cellule $C$ (\cf \cite{Gro} pour la d\'efinition des groupes de cohomologie \`a support dans une sous-vari\'et\'e localement ferm\'ee). 

Comme le faisceau $\call{L}_\lambda$ est $G-$lin\'earis\'e sur $X$, les groupes de cohomologie \`a support $H^i_\call{B}(\call{L}_\lambda) \et H^i_C(\call{L}_\lambda)$ sont des {\it $\goth g-B-$modules}\index{$\goth g-B-$module}, \cad des $\goth g-$modules dont l'action de l'algèbre de Lie de $B$ s'intègre en une action rationnelle du groupe $B$,  (\cf \cite[lem. 11.1]{Kempf}). On peut d\'efinir et calculer leurs multiplicit\'es selon les $\goth g-$modules simples de dimension finie.

En effet,  soient $\mu \in {\mathcal X}$ un caractère dominant, $L(\mu)$ le $\goth g-$ module simple de plus haut poids $\mu$ et $\chi_\mu : Z(\goth g) \to \kk$ son caractère central (l'anneau $Z(\goth g)$ est le centre de l'algèbre enveloppante de $\goth g$). Si $M$ est un $\goth g-$module, on note $M_{\chi_\mu}$ le sous-espace propre g\'en\'eralis\'e associ\'e \`a $\chi_\mu$ (\cf \cite[\S 7.8.15]{Dixmier}). Si $M_{\chi_\mu}$ est de longueur finie, on d\'efinit la multiplicit\'e de $L(\mu)$ selon $M$ par :
$$[M:L(\mu)] := [M_{\chi_\mu} : L(\mu)] \p$$

D'après \cite[pro. 4.6]{toi}, les groupes de cohomologie $H^i_\call{B}(\call{L}_\lambda) \et H^i_C(\call{L}_\lambda)$ ont leur sous-espace propre g\'en\'eralis\'e associ\'e \`a $\chi_\mu$ de longueur finie. De plus, leur multiplicit\'e selon le $\goth g-$module simple $L(\mu)$ est $0$ ou $1$ d'après le lemme \ref{lem:multi0ou1} qui suit.

Rappelons que pour chaque racine sph\'erique $\gamma \in \Sigma_X$, $D_\gamma$ d\'esigne le diviseur limitrophe associ\'e.

Avec ces notations :

\begin{lem}\label{lem:multi0ou1}
 Soient $\mu \in {\mathcal X}$ un caract\`ere dominant et $L(\mu)$ le $\goth g-$module simple de plus haut poids $\mu$. Soient $C$ une cellule de Bialynicki-Birula et $\call{B}$ une $B-$orbite de $X$.

i) Si $i \not= \codim_X (C)$, alors $\ma H^i_{C}(\call{L}_\lambda) = 0$.

Si $i = \codim_X(C)$, alors la multiplicit\'e de $L(\mu)$ dans le $\goth g-$module $\ma H^i_{C}(\call{L}_\lambda)$ est :
$$\left[ H^i_{C}(\call{L}_\lambda) : L(\mu) \right] $$
$$= \left\{\begin{array}{cl} 
1 & \si  \ma w \inv (\mu + \rho) \in \lambda + \rho + \sum_{\gamma \in \Sigma_X \atop C \sub D_\gamma} \ZZ_{> 0} \gamma + \sum_{\gamma \in \Sigma_X \atop C \not\sub D_\gamma} \ZZ_{\le 0} \gamma\\
0 & \sinon 
\end{array}\right.$$
o\`u $w$ est l'\'el\'ement de  $W^Q$ tel que $C= X^+_{w}$.

\saut

ii)  Si $i \not= \codim_X(\call{B})$, alors $\ma H^i_{\call{B}}(\call{L}_\lambda) = 0$.

Si $i = \codim_X(\call{B})$, alors la multiplicit\'e de $L(\mu)$ dans le $\goth g-$module $\ma H^i_{\call{B}}(\call{L}_\lambda)$ est :
$$\left[ H^i_{\call{B}}(\call{L}_\lambda) : L(\mu) \right] $$
$$= \left\{\begin{array}{cl} 
1 & \si \ma w \inv (\mu + \rho) \in \lambda + \rho + \sum_{\gamma \in \Sigma_X \atop \call{B} \sub D_\gamma} \ZZ_{> 0} \gamma + \sum_{\gamma \in \Sigma_X \atop \call{B} \not\sub D_\gamma} \ZZ \gamma\\
0 & \sinon
\end{array}\right.$$
o\`u $w$ est l'\'el\'ement de $W^Q$ tel que $\call{B} \sub  X^+_{w}$.

\saut
\end{lem}

Le point $i)$ d\'ecoule de \cite[th. 4.1 et th 4.4]{toi}. Le point $ii)$ d\'ecoule de \cite[pro. 4.6 et th. 4.4]{toi} et du fait que dans une vari\'et\'e magnifique de rang minimal, tous les points fixes du tore $T$ sont dans l'orbite ferm\'ee $F$ de $G$ (\cf la proposition \ref{pro:pointsfixes}).
\saut

Maintenant, il s'agit de passer des groupes de cohomologie \`a support dans les $B-$orbites aux groupes de cohomologie usuels $$ H^i(X, \call{L}_\lambda) \p$$
C'est l'objet de la prochaine partie :

\section{De la cohomologie \`a support \`a la cohomologie usuelle}

\subsection{Complexe de Grothendieck-Cousin}\label{sec:GC}

Reprenons la filtration 
$$X = Z_0 \con Z_1 \con ...$$
introduite au paragraphe \ref{par:dec} ; on va consid\'erer les groupes de cohomologie \`a support dans les sous-vari\'et\'es localement ferm\'ees  $Z_i \moins Z_{i+1}$ : $H^i_{Z_i \moins Z_{i+1}} (\call{L}_\lambda)$, $i \ge 0$.  Le r\'esultat suivant est d\^u \`a Kempf :
\begin{thm}[{\cite[th. 8.7 (b)]{Kempf}}]

Il existe un complexe de $\goth g-$modules :
\begin{equation}
 ... \to H^i_{Z_i \moins Z_{i+1}}(\call{L}_\lambda) \sta{d^i}{\to} H^{i+1}_{Z_{i+1} \moins Z_{i+2}}(\call{L}_\lambda)\to ...
\end{equation}
dont le  $i-$ème groupe d'homologie est $H^i(X,\call{L}_\lambda)$.

De plus, pour chaque $i$, le $i-$ème terme du complexe se d\'ecompose :
\begin{equation}
H^i_{Z_i \moins Z_{i+1}}(\call{L}_\lambda) = \Plus_\call{B} H^i_\call{B}(\call{L}_\lambda)
\end{equation}
(somme directe sur les $B-$orbites $\call{B}$ de $X$ telles que $\codim_X(\call{B}) = i$) \raisebox{-.2ex}{$\Box$} 
\end{thm}
{\bf Remarque} : {\it Si $\call{B}$ et $\call{B}'$ sont deux $B-$orbites de $X$ de codimensions $i$ et $i+1$, alors :
\begin{liste}
\item la r\'eunion $\call{B} \dij \call{B}'$ est localement ferm\'ee dans $X$,
\item dans $\call{B} \dij \call{B}'$, $\call{B}$ est ouvert et $\call{B}'$ est ferm\'e.
\end{liste}
On en d\'eduit donc une suite exacte longue :
\begin{equation}
... \to H^i_{\call{B} \dij \call{B}'}(\call{L}_\lambda) \to H^i_{\call{B}}(\call{L}_\lambda) \sta{d^i_{\call{B}, \call{B}'}}{\to} H^{i+1}_{\call{B}'}(\call{L}_\lambda) \to ...
\end{equation}
(\cf \cite[pro. 1.9]{Gro}). 
Avec ces notations, dans le th\'eorème ci-dessus, chaque diff\'erentielle 
$$d^i : \Plus_\call{B} H^i_\call{B}(\call{L}_\lambda) \to \Plus_{\call{B}'} H^{i+1}_{\call{B}'}(\call{L}_\lambda)$$
est donn\'e par une \og matrice \fg :
$$\left(d^i_{\call{B}, \call{B}'} \right)_{\call{B}, \call{B}'}$$
o\`u $\call{B}$ et $\call{B}'$ d\'ecrivent, respectivement, les $B-$orbites de codimensions $i$ et $i+1$ dans $X$ (en fait lorsque $\call{B}'$ n'est pas contenue dans l'adh\'erence $\adh{\call{B}}$, $d^i_{\call{B},\call{B}'} = 0$). 
}
\subsection{La partie finie du complexe de Grothendieck-Cousin}

Les $\goth g-$modules $H^d(X,\call{L}_\lambda)$ sont en fait des $G-$modules car le faisceau $\call{L}_\lambda$ est $G-$lin\'earis\'e sur $X$. C'est pourquoi, dans le complexe de Grothendieck-Cousin :\begin{equation}\label{eq:BGG}
{\mathcal GC}^*:... \to \Plus_{\call{B}} H^i_{\call{B}}(\call{L}_\lambda) \to ...
\end{equation}
il suffit de ne pr\^eter attention qu'aux multiplicit\'es selon les $G-$modules simples \ie selon les $\goth g-$modules simples de dimension finie. 

\subsubsection{D\'ecomposition du complexe suivant les caractères centraux}
Comme pr\'ec\'edemment, soient $\mu \in {\mathcal X}$ un caractère dominant, $L(\mu)$ le $G-$module simple de plus haut poids $\mu$ et $\chi_\mu$ son caractère central.

Si $M$ est un  $\goth g-B-$module, alors on note $M_{(\mu)}$ le sous-$T-$espace propre associ\'e au caractère $\mu$ du sous-espace propre g\'en\'eralis\'e $M_{\chi_\mu}$
\index{$M_{(\mu)}$, $\mu-$espace propre de l'espace propre g\'en\'eralis\'e selon le caractère central $\chi_\mu$} \infra{
Par d\'efinition : $M_{(\mu)} = \left\{ m \in \uni_{n >0} \inter_{c \in Z(\goth g)} \ker (c - \chi_\mu(c))^n \tq \qq t \in T, t. m = \mu(t)m \right\}$
}.

Le foncteur ainsi d\'efini : 
$$M \donne M_{(\mu)}$$
v\'erifie :

\begin{pro} 
Si $M' \to M \to M''$ est une suite exacte de $\goth g-B-$modules alors :
\begin{liste}
\item[i)] la suite d'espaces vectoriels 
$$M'_{(\mu)} \to M_{(\mu)} \to M''_{(\mu)}$$
est encore exacte ;

\item[ii)] Si $M_{\chi_\mu}$ est de longueur finie alors l'espace vectoriel $M_{(\mu)}$ a pour dimension :
$$\dim_\kk M_{(\mu)} = [M:L(\mu)] $$
la multiplicit\'e de $M$ selon $L(\mu)$.
\end{liste}
\end{pro}

\begin{dem}
i) cela r\'esulte du fait que $M \donne M_{\chi_\mu}$ est un foncteur exact sur les $\goth g-B-$modules (\cf \cite[\S 7.8.15]{Dixmier})

ii) Les $\goth g-$modules simples qui apparaissent dans une suite de Jordan-H\"{o}lder de $M_{\chi_\mu}$ sont tous de la forme $L(w * \mu)$, $w \in W$, le $\goth g-$module simple de plus haut poids $w * \mu$.

Or, comme le caractère $\mu$ est dominant, parmi les $\goth g-$modules $L(w*\mu)$, seul $L(\mu)$ a une multiplicit\'e non nulle suivant le poids $\mu$ ; en outre cette multiplicit\'e vaut $1$. D'o\`u : $\dim_\kk (M_{(\mu)}) = \dim (M_{\chi_\mu})_{\mu} = [M_{\chi_\mu} : L(\mu)]$. 
\end{dem}

 Par exemple, comme les faisceaux $\call{L}_\lambda$ sont $G-$lin\'earis\'es sur $X$, les groupes de cohomologie \`a support $\ma H^i_\call{B} (\call{L}_\lambda)$ sont bien des $\goth g-B-$modules. En appliquant $M \donne M_{(\mu)}$ au complexe de Grothendieck-Cousin, ${\mathcal GC}^*$,  on obtient un nouveau complexe :\begin{equation}\label{eq:BGGmu}
{\mathcal GC}^*_\mu : ... \to \Plus_{\call{B}} H^i_{\call{B}}(\call{L}_\lambda)_{(\mu)} \to ...
\end{equation}
 qui est plus simple : on va voir en effet que ce nouveau complexe se d\'ecompose en une somme directe de sous-complexes dont on sait calculer l'homologie. 

\subsubsection{D\'ecomposition du complexe suivant les cellules}

Pour qu'un morphisme :
$$d^i_{\call{B}, \call{B'}}(\mu) :  H^i_{\call{B}}(\call{L}_\lambda)_{(\mu)} \to H^{i+1}_{\call{B}'}(\call{L}_\lambda)_{(\mu)}$$

soit non nul, il y a au moins deux conditions n\'ecessaires :
\begin{equation}\label{eq:cond1}
\call{B}' \sub \adh{\call{B}} \et \codim_X(\call{B}') = \codim_X(\call{B})+ 1 = i+1
\end{equation}
car $d^i_{\call{B},\call{B}'} \not=0$ ; et :
\begin{equation}\label{eq:cond2}
\ma H^i_{\call{B}}(\call{L}_\lambda)_{(\mu)} \not = 0 \et H^{i+1}_{\call{B}'}(\call{L}_\lambda)_{(\mu)} \not=0
\end{equation}
 (cela se transforme en une condition combinatoire sur les poids $\lambda$ et $\mu$ d'apr\`es le lemme \ref{lem:multi0ou1}).

En fait, ces deux conditions obligent les deux orbites $\call{B}$ et $\call{B}'$ \`a \^etre dans une m\^eme cellule ; autrement dit :

\begin{lem}\label{lem:memecellule}
Soient $\call{B}$ et $\call{B}'$ deux $B-$orbites de $X$.
Si  $\call{B}$ et $\call{B}'$ sont dans deux cellules de Bialynicki-Birula distinctes \infra{
pour n'importe quelle d\'ecomposition cellulaire donn\'ee par un sous-groupe \`a un paramètre dominant et $X-$r\'egulier.
}, alors l'application :
$$ d^i_{\call{B}, \call{B}'}(\mu) :  H^i_{\call{B}}(\call{L}_\lambda)_{(\mu)} \to H^{i+1}_{\call{B}'}(\call{L}_\lambda)_{(\mu)}$$
est nulle (pour chaque $i$).
\end{lem}

\begin{center}
*
\end{center}

Ce lemme \ref{lem:memecellule} est en fait le lemme clef pour calculer la cohomologie du faisceau $\call{L}_\lambda$ ; admettons-le pour le moment. Sa d\'emonstration, que l'on donne en la section \ref{subsec:demlemcle}, repose sur l'\'etude des courbes irr\'eductibles et $T-$invariantes, dans la vari\'et\'e $X$ (\cf les sections \ref{sec:chaine}, \ref{sec:etu} et \ref{sec:lemcle}).

Un sous-groupe \`a un paramètre $\zeta$, dominant et $X-$r\'egulier, \'etant fix\'e et  :
$$X = \Dij_{x \in X^T} X^+(x)$$\'etant
la d\'ecomposition cellulaire correspondante, pour chaque $i$, soit :
\begin{equation}
Z_i^x := Z_i \cap X^+(x)
\end{equation} 
(c'est la r\'eunion des $B-$orbites $\call{B}$ contenues dans la cellule $X^+(x)$ et telles que $\codim_X(\call{B}) \ge  i$).

On a alors pour tout $i$ une d\'ecomposition :
\begin{equation}
H^i_{Z_i/Z_{i+1}}(\call{L}_\lambda) = \Plus_{x \in X^T} H^i_{Z_i^x/Z_{i+1}^x}(\call{L}_\lambda) \p
\end{equation}

Or, le lemme \ref{lem:memecellule} affirme que pour chaque entier $i$, pour chaque point fixe $x \in X^T$ et pour chaque caractère dominant $\mu \in {\mathcal X}$ :
\begin{equation}
d^i(\mu) \left(  H^i_{Z_i^x/Z_{i+1}^x}(\call{L}_\lambda) _{(\mu)} \right) \sub  H^{i+1}_{Z_{i+1}^x/Z_{i+2}^x}(\call{L}_\lambda) _{(\mu)} 
\end{equation}
 autrement dit, le complexe :
$${\mathcal GC}_\mu^* \; : ...  \to  H^i_{Z_i/Z_{i+1}}(\call{L}_\lambda)_{(\mu)} \sta{d^i(\mu)}{\to}  H^{i+1}_{Z_{i+1}/Z_{i+2}}(\call{L}_\lambda)_{(\mu)} \to ...$$
se d\'ecompose en une somme directe de complexes :
\begin{equation}\label{eq:mux}
{\mathcal GC}^*_\mu = \Plus_{x \in X^T} {\mathcal GC}_{\mu,x}^*
\end{equation}
o\`u pour chaque $x \in X^T$, on d\'efinit le complexe :
$${\mathcal GC}_{\mu,x}^* \; : ...  \to  H^i_{Z_i^x/Z_{i+1}^x}(\call{L}_\lambda)_{(\mu)} \sta{d^i(\mu)}{\to}  H^{i+1}_{Z_{i+1}^x/Z_{i+2}^x}(\call{L}_\lambda)_{(\mu)} \to ...$$

Mais, c'est aussi, pour chaque $x \in X^T$, la composante suivant le caractère dominant $\mu$ du complexe :
$${\mathcal GC}^*_x \; : ...  \to H^i_{Z_i^x/Z_{i+1}^x}(\call{L}_\lambda) \sta{d^i_x}{\to} H^{i+1}_{Z_{i+1}^x/Z_{i+2}^x}(\call{L}_\lambda)  \to ...$$
correspondant \`a la filtration : 
$$X^+(x) = Z_0^x \con Z_1^x \con ...$$

{\bf Remarque :} {\it Chaque complexe ${\mathcal GC^*}_x$ a donc pour homologie :
\begin{equation}\label{eq:GCc}
H^* ({\mathcal GC^*}_x) = H^*_{X^+(x)}(\call{L}_\lambda)
\end{equation}
et ces groupes de cohomologie \`a support sont donn\'es par le lemme \ref{lem:multi0ou1}}.
\subsection{Calcul de la cohomologie  usuelle}

\'Etant donn\'e que les groupes de cohomologie $\ma H^d(X,\call{L}_\lambda)$ sont non seulement des $\goth g-$modules mais aussi des $G-$modules, il r\'esulte de (\ref{eq:mux}) et  (\ref{eq:GCc}), en tout degr\'e $d$ et pour tout poids dominant $\mu$, l'\'egalit\'e :
\begin{equation}\label{eq:supssup}
 \left[ H^d(X,\call{L}_\lambda): L(\mu) \right] =
\sum_{w \in W^Q }\left [  H^d_{X^+_w}(\call{L}_\lambda) : L(\mu) \right ] \p
\end{equation}

Par cons\'equent, pour terminer la d\'emonstration du th\'eorème \ref{thm:principal}, il s'agit maintenant de calculer la multiplicit\'e selon le $\goth g-$module simple $L(\mu)$ des $\goth g-B-$modules 
$$H^d_{X^+_w}(\call{L}_\lambda) \p$$

Fixons un poids $\mu$ dominant. Rappelons que $\rho$ d\'esigne la demi-somme des racines positives.

D'après le lemme \ref{lem:multi0ou1}, si on choisit $w \in W^Q$, \cad un point fixe $x = w\zzz \in X^T$, alors, d'une part :
$$
\left[ H^d_{X^+_w}(\call{L}_\lambda) : L(\mu) \right] = 1$$
si les deux conditions suivantes sont v\'erifi\'ees :
\begin{equation}\label{eq:codim}
\codim_X(X^+_w) = d
\end{equation} \et 
\begin{equation}\label{eq:combi}
w\inv (\mu + \rho)  \in \lambda + \rho + \sum_{\gamma \in \Sigma_X \atop X^+_w \sub D_\gamma} \ZZ_{> 0} \gamma + \sum_{\gamma \in \Sigma_X \atop X^+_w \not\sub D_\gamma} \ZZ_{\le 0} \gamma 
\end{equation}
et d'autre part, dans tous les autres cas :

$$\left[ H^d_{X^+_w}(\call{L}_\lambda) : L(\mu) \right] = 0 \p$$
\begin{center}
*
\end{center}
Au passage, notons que la condition (\ref{eq:codim}) entra\^ine que $w = w_\nu$ avec :
$$\nu := w\inv(\mu+\rho)-\rho \in \pic(X) \p$$

A priori, suivant la proposition \ref{pro:coce}, la codimension de la cellule $X^+_w$ et l'inclusion (ou non) $X^+_w \sub D_\gamma$ d\'ependent du sous-groupe \`a un paramètre dominant et $X-$r\'egulier $\zeta$ choisi pour d\'efinir la d\'ecomposition cellulaire. Mais, en fait, ce n'est pas le cas lorque $w$ est de la forme $w_\nu$ pour un poids $\nu \in \pic(X)$.

En effet :

\begin{pro}\label{pro:wparti}
Pour toute racine sph\'erique $\gamma$ et pour tout $\nu \in \pic(X)$, on a :
$$\cg w_\nu \gamma , \zeta \cd > 0 \equi (\nu + \rho , \gamma) > 0 \p$$
\end{pro}

\begin{dem}
 Nous verrons plus loin qu'il existe deux racines positives orthogonales  $\alpha$ et $\beta$ telles que $\alpha + \beta \in \Z_{>0} \gamma$ et $\cg \nu + \rho , \alpha^\ch\cd = \cg \nu + \rho , \beta^\ch \cd$ pour tout $\nu \in \pic(X)$ (\cf la remarque qui suit le lemme \ref{lem:sg} page \pageref{lem:sg}). D'un c\^ot\'e, on en d\'eduit que :
$$\cg w_\nu \gamma, \zeta \cd \et  \cg w_\nu(\alpha) , \zeta \cd + \cg w_\nu(\beta) , \zeta \cd$$
sont de m\^eme signe. 

D'un autre c\^ot\'e, comme le sous-groupe \`a un paramètre $\zeta$ est dominant et $X-$r\'egulier, on a les \'equivalences :
\begin{equation}\label{eq:a1}
\cg  w_\nu(\alpha) , \zeta \cd > 0 \equi w_\nu(\alpha) > 0
\end{equation}

\begin{equation}\label{eq:b}
  \cg w_\nu(\beta) , \zeta \cd > 0 \equi w_\nu(\beta) > 0 \p
\end{equation}
Mais par d\'efinition de $w_\nu$, $w_\nu(\nu + \rho) \in {\mathcal X} $ est un poids dominant r\'egulier et donc :
$$w_\nu(\alpha) > 0 \equi (w_\nu(\nu + \rho) , w_\nu(\alpha)) > 0$$
\begin{equation}\label{eq:ar}
\equi (\nu + \rho, \alpha) > 0 \p
\end{equation}
De m\^eme :\begin{equation}\label{eq:br}
w_\nu(\beta) > 0 \equi(\nu + \rho, \beta) > 0 \p
\end{equation}

Puisque $\cg \nu + \rho , \alpha^\ch \cd = \cg \nu + \rho , \beta^\ch \cd$, on a : $w_\nu (\alpha) > 0 \equi w_\nu(\beta) > 0$. Il r\'esulte alors de (\ref{eq:a1}), (\ref{eq:b}), (\ref{eq:ar}), (\ref{eq:br}) que :
$$\cg w_\nu \gamma , \zeta \cd > 0 \equi  \cg w_\nu (\alpha) , \zeta \cd + \cg w_\nu(\beta) , \zeta \cd >0$$
$$ \equi (\nu + \rho , \alpha) + (\nu + \rho , \beta) > 0 \equi (\nu + \rho , \gamma) > 0 \p$$

\end{dem} 

\begin{center}
*
\end{center}

Un poids dominant $\mu$ et un \'el\'ement $w$ de $W^Q$ \'etant toujours fix\'es, on garde la notation $\nu : = w \inv (\mu + \rho) -\rho$. 

D'une part $\nu + \rho$ est un caractère r\'egulier et d'autre part, si les conditions (\ref{eq:codim}) et (\ref{eq:combi}) sont remplies alors, d'après la proposition \ref{pro:wparti} et avec les notations de l'\'enonc\'e du th\'eorème principal \ref{thm:principal}, on trouve : 

$$\nu \in (\lambda + R_J) \cap \Omega_J \et l(\nu) + |J| = d$$
avec $J := \{\gamma \in \Sigma_X \tq (\nu + \rho , \gamma) < 0\}$.

R\'eciproquement, si $\nu \in {\mathcal X}$ (tel que $\nu^+$ existe, \ie tel que $\nu + \rho$ soit r\'egulier) est dans un ensemble de poids :
$$(\lambda + R_J) \cap \Omega_J$$
pour une certaine partie $J$ de $\Sigma_X$, alors le $\goth g-$module :
$$H^d_{X^+_w}(\call{L}_\lambda)$$
avec $d := l(\nu) + |J|$ et  $w:=w_\nu \in W^Q$ a pour multiplicit\'e $1$ selon le $\goth g-$module simple $L(\nu^+)$.

La formule du th\'eorème principal \ref{thm:principal} :

$$H^d(X,\call{L}_\lambda) \iso \Plus_{J \sub \Sigma_X} \Plus_{\mu \in (\lambda + R_J) \cap \Omega_J \atop{\mu + \rho \mbox{ \scriptsize r\'{e}gulier } \atop
l(\mu) +  |J|= d} } L(\mu^+) $$
d\'ecoule alors de (\ref{eq:supssup}), au d\'ebut de cette section.

Le th\'eorème \ref{thm:principal} est presque d\'emontr\'e.

\vskip .5cm

Il ne reste plus, en effet qu'\`a v\'erifier le lemme clef \ref{lem:memecellule}. 

Pour cela, on va, dans les deux sections suivantes, d'une part, relier certains points fixes du tore $T$ par des cha\^ines de courbes $T-$invariantes contenues dans la vari\'et\'e $X$ ; et d'autre part \'etudier ces courbes $T-$invariantes.

Cela permettra de d\'emontrer le lemme \ref{lem:memecellule}, dans la section \ref{sec:lemcle}.


\section{Cha\^ines de courbes qui relient deux points fixes du tore}\label{sec:chaine}

Dans cette section \ref{sec:chaine}, on suppose seulement que $X$ est une vari\'et\'e projective munie d'une action \og lin\'eaire \fg\ du tore $T$ au sens suivant : il existe un $T-$module $V$ tel que $X \sub \PP(V)$ et tel que l'action de $T$ sur $X$ provienne de celle de $T$ sur $\PP(V)$.

 On rappelle qu'un sous-groupe \`a un param\`etre $\nu : \kk^* \to T$
est  $X-$r\'egulier si les ensembles de points fixes $$X^\nu =\{x \in X \tq
\qq a \in k^* , \; \nu(a) . x = x \}$$ et 
$$X^T= \{x \in X \tq \qq t \in T, \;
t.x = x\}$$ sont les m\^emes.

\subsection{Orientation des courbes}\label{sec:odc}

Soit $\nu$ un sous-groupe \`a un param\`etre de $T$ qui est $X-$r\'egulier.

\begin{defin}
 Si
$x,x' \in X^T$ sont deux points fixes de $T$, on notera $$x \sta{c}{\to} x'\index{$x \sta{c}{\to} x'$, courbe $c$ reliant $x$ \`a $x'$} \mbox{ ou simplement } x \to x'$$ si la
vari\'et\'e $X$ contient une courbe $c$ ferm\'ee, irr\'eductible,
{\it stable par $T$}, \og allant de $x$ \`a $x'$ \fg\  \ie : il existe $y \in c$ tel que :
$$\lim_{a \to 0}
\nu(a).y = x \et \lim_{a \to \infi} \nu(a).y = x' \p$$

Dans cette situation, on notera :
$$c(0) := x \et c(\infi) :=x' \p$$
\end{defin}

Ces notations d\'ependent du sous-groupe \`a un paramètre $X-$r\'egulier, $\nu$, choisi : si on remplace $\nu$ par $-\nu$, on \'echange $c(0)$ et $c(\infi)$.

\vskip 2em

{\bf Exemple :} Si on prend pour $X$ la vari\'et\'e de drapeaux $G/B^-$, les
points fixes du tore $T$ sont les $wB^-/B^-$ o\`u $w \in W$. Pour un
sous-groupe \`a un param\`etre $\nu$, {\it dominant} et $X-$r\'egulier (\ie r\'egulier), et pour $w,w' \in W$, on a l'\'equivalence :
$$wB^-/B^- \to w'B^-/B^- $$
$$\equi w' =  w s_\alpha$$
pour une racine
positive $\alpha$ telle que $w(\alpha) >0$ ($s_\alpha \in W$ d\'esigne la r\'eflexion associ\'ee \`a la racine $\alpha$). En particulier,  si $wB^-/B^- \to w'B^-/B^-$ alors $w' > w$ pour l'ordre de Bruhat sur $W$.

En effet, les courbes irr\'eductibles et $T-$invariantes de $G/B^-$ qui passent par $wB^-/B^-$ sont de la forme $\adh{wU_\alpha B^-/B^-}$ o\`u $\alpha$ est une racine positive et $U_\alpha$ le sous-groupe unipotent de dimension $1$ associ\'e. Or, les $T-$points fixes de la courbe $\adh{wU_\alpha B^-/B^-}$ sont $wB^-/B^-$ et $ ws_\alpha B^-/B^-$. De plus, si $w(\alpha) > 0$, alors $wB^-/B^- \to ws_\alpha B^-/B^-$, tandis que si $w(\alpha) <0$, alors $ws_\alpha B^-/B^- \to wB^-/B^-$.

\subsection{Liaison entre deux points fixes du tore}

Avec les notations de la section pr\'ecedente, on a le r\'esultat suivant :

\begin{lem}\label{lem:chemin} Soient $V$ un $T-$module de dimension
finie et $X$ une sous-vari\'et\'e ferm\'ee de l'espace projectif
$\PP(V)$. On munit $X$ de l'action de $T$ induite et on suppose que
$X$ n'a qu'un nombre fini de points fixes pour cette action.

On fixe un sous-groupe \`a un paramètre $X-$r\'egulier, $\nu$.  
Soient $x', x$ deux points fixes de $X$.

Si $ x'$ est dans l'adh\'erence de la cellule $X^+(x)$ alors il existe des points fixes $x_0,...,x_N$ dans $X$ tels que :

$$x=x_0 \to x_1 \to ... \to x_N = x' \p$$
\end{lem}

{\bf Remarque :} Si $X$ est la vari\'et\'e de drapeaux $G/B^-$ et si on choisit un sous-groupe \`a un paramètre $\nu$, dominant et r\'egulier, alors :
$$w'B^-/B^- \in \adh{X^+_w} \equi w \le w' \mbox{ pour l'ordre de Bruhat}$$
$$\impliq \exist \alpha_1,..., \alpha_n \in \Phi^+ \,, w < ws_{\alpha_1} <...< ws_{\alpha_1}...s_{\alpha_n} = w'$$
et on a bien une cha\^ine de courbes :
$$wB^-/B^- \to ... \to w'B^-/B^- \p$$

\begin{center}
*
\end{center}

\begin{dem}[du lemme]
Toute composante irr\'eductible de $\adh{X^+(x)}$ est l'adh\'erence d'une composante irr\'eductible de $X^+(x)$. En particuler, toutes les composantes irr\'eductibles de $\adh{X^+(x)}$ contiennent $x$. On peut donc supposer que $X = \adh{X^+(x)}$, quitte \`a remplacer $X$ par une composante irr\'eductible de $\adh{X^+(x)}$ qui contient $x'$.

Si $z$ est un point fixe de $X$, on pose :

$$X^-(z):= \{ y \in X \tq \lim_{a \to \infi} \nu(a) . y = z \} \p$$

\'Etant donn\'es deux points fixes distincts $z,z'$, on \'ecrira :
$$z
\gamb z' \si  X^+(z) \cap X^-(z') \not= \vide \p$$

Via le sous-groupe \`a un paramètre $\nu$, le groupe $\kk^*$ agit sur l'espace vectoriel $V$ avec des poids entiers que l'on peut ordonner. 

Pour un point fixe $x$ de $T$, notons $\nu_x$ le poids d'un vecteur propre $v \in
V$ tel que $x = [v]$ (la classe du vecteur $v$ dans l'espace projectif). 

On s'aperçoit que si $x \gamb x'$, alors
$\nu_{x} < \nu_{x'}$ (il suffit de le v\'erifier dans $\PP(V)$). 

On en d\'eduit que dans la sous-vari\'et\'e $X$ de
$\PP(V)$, il n'y a pas de  \og boucles \fg\ \ie : s'il existe une cha\^ine $x_0
\gamb x_1 \gamb ... \gamb x_N = x_0$, alors $N=0$. 

En particulier,
puisque $X$ n'a qu'un nombre fini de points fixes, il existe  une
cha\^ine aboutissant \`a $x'$ : \begin{equation}\label{eq:aboumax}
 x_0 \gamb ... \gamb x_N=x'
\end{equation}  avec un
entier $N$ maximal.

On va montrer que $x_0 = x$ : 

si $y \in X^-(x_0) \moins \{x_0\}$,  alors
$x_0 \not= \lim_{a \to 0} \nu(a).y \gamb x_0$ ce qui contredit la maximalit\'e de
$N$. Ainsi, $X^-(x_0) = \{x_0\}$. Or, d'apr\`es \cite{Konarski}[th. 3 et
bas de la p. 299] : $$ \dim X^-(x_0) + \dim X^+(x_0) \ge \dim X \p$$
Il s'ensuit que $\dim X^+(x_0) = \dim X$ et donc que $\adh{X^+(x_0)} =
X$. Comme $X^+(x_0)$ et $X^+(x)$ sont ouverts, cela
entra\^ine que $X^+(x_0)$ rencontre $X^+(x)$ et par cons\'equent que $x_0 = x$.

Montrons maintenant que dans la cha\^ine maximale (\ref{eq:aboumax}), $x_i \to x_{i+1}$ pour tout $0 \le i \le
N-1$. La courbe $\adh{\kk^* . y}$ pour un $y \in X^+(x_i) \cap X^-(x_{i+1})$ est bien une courbe allant de $x_i$ \`a $x_{i+1}$ mais a
priori elle n'est que $\kk^*-$invariante. 

Cependant, si $z$ est un point fixe du tore $T$ dans $\adh{X^+(x_{i}) \cap X^-(x_{i+1})}$ alors, comme $z \in \adh{X^+(x_i)}$, il existe une cha\^ine partant de $x_i$ et aboutissant \`a $z$ :
$$x_{i} \gamb ... \gamb z$$
et, de m\^eme, comme $z\in \adh{X^-(x_{i+1})}$, il existe une cha\^ine aboutissant \`a $x_{i+1}$ et partant de $z$ :
$$z \gamb ... \gamb x_{i+1} \p$$

Par maximalit\'e de la cha\^ine $x_0 \gamb ... \gamb x_i \gamb x_{i+1} \gamb ...\gamb x_N$, il est n\'ecessaire que $z = x_i$ ou $x_{i+1}$. En particulier, il n'y a que deux points fixes dans $\adh{X^+(x_{i}) \cap X^-(x_{i+1})}$. D'apr\`es \cite{Ros} ou \cite[cor. 1 p. 497]{BiBi}, la sous-vari\'et\'e ferm\'ee $\adh{X^+(x_{i+1}) \cap X^-(x_i)}$ est donc irr\'eductible et de dimension $1$. C'est donc une courbe ferm\'ee allant de $x_i$ \`a $x_{i+1}$ qui, de plus, est $T-$invariante ; ainsi $x_i \to x_{i+1}$.
\end{dem}

On revient maintenant au cas o\`u la vari\'et\'e $X$ est magnifique et de rang minimal.

\section{\'Etude de quelques courbes invariantes par le tore}\label{sec:etu}
Dans la vari\'et\'e magnifique de rang minimal, $X$, il n'y a, en fait, qu'un nombre fini de courbes $T-$invariantes (\cf le lemme \ref{lem:clambda} ci-dessous). Parmi ces courbes, celles qui sont incluses dans la $G-$orbite ferm\'ee $F$ sont des courbes $T-$invari\-antes d'une vari\'et\'e de drapeaux (courbes bien connues : \cf \cite{Springer97}, \cf aussi l'exemple de la section \ref{sec:odc}). On va plut\^ot s'int\'eresser aux courbes $T-$in\-variantes de $X$ qui sortent de l'orbite ferm\'ee $F$ ; on associera de telles courbes aux racines sph\'eriques $\gamma \in \Sigma_X$. Après avoir donn\'e des exemples en rang $1$, on \'etablira quelques propri\'et\'es des racines sph\'eriques, importantes pour d\'emontrer le lemme clef  \ref{lem:memecellule}.

\subsection{Finitude des courbes invariantes par le tore}

La vari\'et\'e magnifique $X$ \'etant de rang minimal, ses racines sph\'eriques ne sont pas des multiples de racines (\cf le lemme \ref{lem:racnonsph} (section \ref{sec:rac-sph})). En particulier, les poids de l'espace tangent $T_{\zzz} X$ sont deux \`a deux non proportionnels (en effet, ces poids sont les racines sph\'eriques $\gamma \in \Sigma_X$ et certaines racines positives $\alpha \in \Phi$).

Par cons\'equent, il n'y a qu'un nombre fini de courbes irr\'eductibles et $T-$invariantes passant par $\zzz$ ; plus pr\'ecis\'ement :

\begin{lem}\label{lem:clambda}
Soit $\lambda$ un poids de l'espace tangent $T_{\zzz} X$. Il existe une unique courbe $C_\lambda$ dans $X$ irr\'eductible et $T-$invariante telle que :
$$\zzz \in C_\lambda \et T_{\zzz}C_\lambda = (T_{\zzz}X)_\lambda$$
l'espace propre de $T_{\zzz}X$ de poids $\lambda$. 

De plus, la courbe $C_\lambda$ est la composante connexe de $X^{(\ker \lambda)^\circ}$ qui contient $\zzz$. Enfin, toute courbe $T-$invariante $C$ de $X$ qui contient $\zzz$ est l'une des courbes $C_\lambda$, en particulier est lisse.
\end{lem}

{\bf Remarques :} 

\begin{liste}
\item[i)] On a donc une bijection entre les poids de l'espace tangent $T_{\zzz} X$ et les courbes irr\'eductibles et $T-$invariantes de $X$ passant par $\zzz$.

\item[ii)] Comme tous les points fixes du tore $T$ sont dans la $G-$orbite ferm\'ee $F$, toute courbe irr\'eductible et $T-$invariante de $X$ est conjugu\'ee par un $w \in W^Q$ \`a une courbe irr\'eductible et $T-$invariante passant par $\zzz$.

\item[iii)] Dans le cadre des vari\'et\'es de drapeaux, ce r\'esultat est bien connu (\cf par exemple \cite{Springer97}) : si $X$ est la vari\'et\'e de drapeaux $G/B^-$, alors les courbes $T-$invariantes de $X$ passant par $\zzz = B^-/B^-$ sont les adh\'erences $\adh{U_\alpha B^- / B^-}$, o\`u $\alpha$ d\'ecrit l'ensemble des racines positives ($U_\alpha$ d\'esigne le sous-groupe unipotent de $G$ associ\'e \`a la racine $\alpha$). En effet, les poids de l'espace tangent $T_{\zzz} G/B^-$ sont exactement les racines positives et, pour toute racine positive $\alpha$, la vari\'et\'e $U_\alpha B^- / B^-$ est de dimension $1$ et contenue dans $(G/ B^-)^{(\ker \alpha)^\circ}$.
\end{liste}

\begin{dem}
Soit $S := (\ker \lambda)^\circ$ ; c'est un sous-tore de $T$ de codimension $1$.

{\bf existence :} D'après \cite{BiBi}, la sous-vari\'et\'e ferm\'ee $X^S$ des points fixes de $S$ est lisse et pour tout $x \in X^S$ :
$$T_x(X^S) = (T_xX)^S= \left\{ v \in T_x X \tq \qq s \in S ,\; s. v = v \right\} \p$$

Donc la composante connexe $C$ de $X^S$ qui contient $\zzz$ est une vari\'et\'e lisse, irr\'eductible et $T-$invariante v\'erifiant :
$$T_{\zzz} C = T_{\zzz} X^S = (T_{\zzz} X)^S = \Plus_{\mu \in {\mathcal X}} (T_{\zzz} X)_\mu^S$$
mais, comme $S = (\ker \lambda )^\circ$, pour un caractère $\mu \in {\mathcal X}$ :
$$(T_{\zzz} X)_\mu^S = \{0\}$$ si $\mu$ n'est pas proportionnel \`a $\lambda$ ; comme de plus, $\lambda$ est le seul poids de $T_{\zzz}X$ proportionnel \`a $\lambda$, on trouve :
$$  T_{\zzz} C = (T_{\zzz} X)_\lambda$$
qui est de dimension $1$. Par cons\'equent, $C$ est bien une courbe.

\vskip .3cm

{\bf unicit\'e :} Si $C$ est une courbe irr\'eductible et $T-$invariante qui contient $\zzz$ et si $T_{\zzz}C = (T_{\zzz}X)_\lambda$, alors :
$$(T_{\zzz}C^S) = (T_{\zzz} C)^S = (T_{\zzz}X)_\lambda^S$$
$$=(T_{\zzz} X)_\lambda = T_{\zzz} C \p$$
Mais alors, $C^S = C$ et $C \sub X^S$ : \cad $C$ est la composante connexe de $X^S$ passant par $\zzz$. 

En fait, si $C$ est une courbe irr\'eductible et $T-$invariante de $X$, alors $C$ est fix\'ee point par point par un sous-tore $S$ de $T$ de codimension $1$. Donc si la courbe $C$ passe par le point $\zzz$, on a l'inclusion d'espaces tangents :
$$T_{\zzz} C \sub T_{\zzz} X^S = (T_{\zzz} X)^S \p$$
Comme les poids de $T$ dans $T_{\zzz} X$ sont deux \`a deux non proportionnels, $T_{\zzz} C$ est forc\'ement de dimension $1$ et en cons\'equence, $C$ est lisse. 
\end{dem}

\subsection{Courbes et r\'eflexions associ\'ees aux racines sph\'eriques}\label{subsec:cgamma}

Gr\^ace au lemme \ref{lem:clambda}, on peut poser la d\'efinition suivante :

\begin{defin}\label{def:cgamma}
Pour toute  racine sph\'erique $\gamma \in \Sigma_X$, soit $\ma
C_\gamma$\index{$C_\gamma$, courbe associ\'ee \`a la racine sph\'erique $\gamma$} l'unique courbe irr\'eductible, $T-$invariante, contenant $\zzz$, contenue dans $X$ telle que :
$$T_{\zzz} C_\gamma = (T_{\zzz}X)_\gamma \p$$
\end{defin}

\begin{center}
*
\end{center}

D'après le lemme \ref{lem:clambda}, les courbes $C_\gamma$ sont exactement les courbes irr\'eductibles, $T-$invariantes, passant par $\zzz$ et qui sortent de $F$.

D'un autre c\^ot\'e, chaque courbe $C_\gamma$ est lisse, isomorphe \`a $\PP^1$ et contient deux points fixes du tore. Cette remarque permet d'associer un \'el\'ement de $W$ \`a chaque racine sph\'erique :

\begin{defin}[l'\'el\'ement $s_\gamma$]\label{def:sg}
Pour chaque racine sph\'erique $\gamma \in \Sigma_X$, soit $s_\gamma\index{$s_\gamma$, r\'eflexion sph\'erique}$ l'unique \'el\'ement de $W^Q$ tel que :
$$\zzz \et s_\gamma \zzz$$
sont les $2$ points fixes du tore dans la courbe $C_\gamma$. On appellera $s_\gamma$ une r\'eflexion sph\'erique\index{r\'eflexion sph\'erique}.
\end{defin}

{\bf Remarque :} 
En fait, l'\'el\'ement $s_\gamma$ n'est pas une r\'eflexion dans $W$ ; c'est seulement une r\'eflexion comme transformation du r\'eseau $\pic(X)$. Effectivement : nous verrons plus tard (\cf le lemme \ref{lem:sg}) que $s_\gamma^2 = 1$ et d'un autre c\^ot\'e :
\begin{equation}\label{eq:sgl}
s_\gamma \lambda - \lambda \in \Z \gamma
\end{equation}
pour tout $\lambda \in \pic(X)$.

En effet, d'après \cite[pro. p. 377]{B:Mori}, la diff\'erence des poids des fibres en $\zzz$ et $s_\gamma \zzz$ du faisceau inversible $\call{L}_\lambda$ est un multiple de $\gamma$, le poids associ\'e \`a la courbe $C_\gamma$ reliant $\zzz$ \`a $s_\gamma \zzz$. Or, les droites
$$\call{L}_\lambda\res{\zzz}  \et \call{L}_\lambda\res{s_\gamma \zzz}$$
ont pour poids respectifs $\lambda$ et $s_\gamma \lambda$, d'o\`u (\ref{eq:sgl}).

\subsection{Exemples en rang  un}\label{ssec:exrun}

On utilisera certaines propri\'et\'es des racines sph\'eriques $\gamma$ et des \'el\'ements $s_\gamma$ associ\'es (\cf le lemme \ref{lem:sg} de la section \ref{sec:sg}).  Pour \'etablir ces propri\'et\'es, on se ramènera au cas o\`u $X$ est une vari\'et\'e magnifique de rang $1$, et m\^eme une des vari\'et\'es magnifiques \'etudi\'ees en exemple ci-dessous. 

\subsubsection{Diagrammes}
Parmi les vari\'et\'es magnifiques de rang $1$ de \cite[table 2 p. 67]{Akh}, celles qui sont de rang minimal correspondent, avec les notations de \cite{Akh}, aux diagrammes :

$$ D_n \;:\;\dnn \;\;\;\; (n\ge 3) \mbox{ ou } D_2 \; :\; \axaa$$
$$ B_3\;:\;\spii$$
Elles sont munies respectivement de l'action  de $Spin_{2n}(\kk)$ et de $Spin_7(\kk)$.
\saut

{\bf Remarque :} Les \og ronds blancs \fg\ des diagrammes correspondent aux racines simples du groupe r\'eductif $Q/R_u(Q)$ ($R_u(Q)$ \index{$R_u(Q)$, radical unipotent de $Q$} \'etant le radical unipotent de $Q$). En particulier, les ensembles $W^Q$ sont respectivement :

$$ \{w \in W \tq \qq i \ge 2 , w(\alpha_i) > 0 \} \mbox{ et } W \;\;;$$
$$ \{w \in W \tq w(\alpha_1) > 0, w(\alpha_2) > 0\} \p$$

\subsubsection{Description de l'action}

Les vari\'et\'es magnifiques associ\'ees aux diagrammes ci-dessus sont donn\'ees dans \cite[pp. 66-67]{Akh}.  Dans le cas du diagramme $D_n$, il s'agit de l'espace projectif $\PP^{2n-1}$ ($n \ge 2)$ et de la quadrique de dimension $2n-1$ : $$ \qu_{2n-1} = \left\{ [z] \in \PP^{2n} \tq z_0^2 = \sum_{i=1}^n z_i z_{n+i} \right\}$$ munies d'une action du groupe $Spin_{2n}(\kk)$.

Dans le cas du diagramme $B_3$, il s'agit de l'espace projectif $\PP^7$ et de la quadrique $\qu_7$ (\cf ci-dessus) munies d'une action du groupe $Spin_7(\kk)$.

Pour d\'ecrire les actions, on note $J_{2n}$\index{$J_{2n}$} la matrice carr\'ee d'ordre $2n$ suivante :

$$\left( \begin{array}{ccc}
{\mbox{\Huge $0$}} & & 1\\
 & \antidiago & \\
1 & &  {\mbox{\Huge $0$}}
\end{array}
\right)$$

et $SO(J_{2n})$ le groupe :
$$\left\{ g \in SL(2n,\kk) \tq {}^tgJ_{2n}g = J_{2n}\right\}  $$
qui est isomorphe \`a $SO(2n,\kk)$.
Le groupe $ SO(J_{2n})$ agit sur l'espace projectif $\PP^{2n-1}$ via l'action naturelle de $SL(2n,\kk)$ et sur la quadrique $\qu_{2n-1}$ par :
$$\qq g \in SO(J_{2n}), \; \qq [z_0:z] \in \PP(\kk \croi \kk^{2n}), \; g.[z_0:z] = [z_0:g.z] \p$$

Comme $Spin(2n,\kk)$ est le rev\^etement universel de $SO(J_{2n}) \iso SO(2n,\kk)$ et comme $Spin(7,\kk) \sub SO(J_8)$ (gr\^ace \`a la repr\'esentation spinorielle de dimension $8$ de $Spin(7,\kk)$), on en d\'eduit les actions de $Spin(2n,\kk) \et Spin(7,\kk)$.
\subsubsection{Orbites}

Les vari\'et\'es magnifiques $\PP^{2n-1}$ et $\qu_{2n-1}$, de rang minimal $1$, ont chacune exactement $2$ orbites  de $G = Spin_{2n}(\kk)$ (ou $Spin(7,\kk)$ pour $\PP^7$ et $\qu_8$).

Dans le cas de $\PP^{2n-1}$, l'orbite ouverte et l'orbite ferm\'ee sont respectivement :
$$\call{O} = \left\{[z_1:...:z_{2n}] \in \PP^{2n-1} \tq \sum_{i=1}^n z_i z_{n+i} \not= 0 \right\} $$
$$\et F=\left\{[z_1:...:z_{2n}] \in \PP^{2n-1} \tq \sum_{i=1}^n z_i z_{n+i} = 0 \right\} $$
dans le cas de $\qu_{2n-1}$ ce sont respectivement : 
$$\call{O} = \left\{[z_0 : z_1:...:z_{2n}] \in \qu_{2n-1}\tq z_0 \not= 0 \right\} $$
$$\et F= \left\{[z_0 : z_1:...:z_{2n}] \in \qu_{2n-1}\tq z_0 = 0 \right\} \p$$

Notons $r : Spin_{2n}(\kk) \to SO(J_2n) $ le rev\^etement universel.

Pour le type $D_n$, on prend pour tore maximal le tore :
$$T:= r\inv\left( T_{2n} \cap SO(J_{2n}) \right)$$
o\`u $T_{2n}$ est le sous-groupe des matrices diagonales de $SL(2n,\kk)$. Pour sous-groupe de Borel $B^-$, on prend le groupe :
$$r\inv(B^-_{2n} \cap SO(J_{2n}))$$
o\`u $B^-_{2n}$ est le sous-groupe des matrices triangulaires inf\'erieures de $SL(2n,\kk)$.

Pour le type $B_3$ il suffit de remplacer ci-dessus $r$ par l'inclusion $i : Spin(7,\kk) \to SO(J_8)$.

\subsubsection{Racines et r\'eflexions sph\'eriques}

Comme les vari\'et\'es $X = \PP^{2n-1}$ ou $\qu_{2n-1}$ sont des $G-$vari\'et\'es magnififiques de rang minimal $1$ (pour $G = Spin_{2n}(\kk) \ou Spin_{7}(\kk)$ (si $n=3$)), elles n'ont qu'une racine sph\'erique, $\gamma$. 

En particulier, par l'unique point fixe de $B^-$ dans $X$, not\'e $\zzz$, il ne passe qu'une seule courbe $C_\gamma$, irr\'eductible, $T-$invariante et non contenue dans l'orbite ferm\'ee $F$.

On explicite $C_\gamma$ dans le cas de $\PP^{2n-1}$ :
$$C_\gamma=\left\{ [x:0:...:0:y] \in \PP^{2n-1}\right\}$$ et dans le cas de $\qu_{2n-1}$ :
$$C_\gamma = \left\{[z:x:0:...:0:y] \tq z^2 = xy \right\} \p$$

La courbe $C_\gamma$ joint le point $\zzz$ au point $s_\gamma \zzz$ (ce qui d\'efinit l'\'el\'ement $s_\gamma\in W^Q$). 

Dans le tableau ci-après, on rassemble pour chaque vari\'et\'e magnifique, les points $\zzz$ et $s_\gamma \zzz$.

$$
\begin{array}{|c|c|c|c|c|}
\hline
\mbox{ diagrammes } &\multicolumn{2}{|c|}{\scriptstyle \dnn \;\;\;\; (n\ge 3) \mbox{ ou } \axaa} & \multicolumn{2}{|c|}{\scriptstyle \spii}\\
\hline
G & \multicolumn{2}{|c|}{Spin_{2n}(\kk)} & \multicolumn{2}{|c|}{Spin_7(\kk)}\\
\hline
X & \qu_{2n-1} & \PP^{2n-1} & \qu_7 & \PP^7 \\
\hline

\zzz \;, \atop s_\gamma \zzz&\scriptstyle [0:...:0:1],\atop [0:1:0:...:0] &\scriptstyle [0:...:0:1],\atop [1:0:...:0] &\scriptstyle [0:...:0:1], \atop [0:1:0:...:0]&\scriptstyle
[0:...:0:1],\atop [1:0:...:0]\\
\hline

\end{array}
$$

\subsection{Propri\'et\'es des r\'eflexions sph\'eriques}\label{sec:sg}

Dans le lemme suivant, sont rassembl\'ees les propri\'et\'es des $s_\gamma$ (\cf la d\'efinition \ref{def:sg}) qui vont servir plus loin :

\begin{lem}\label{lem:sg}
Si $\gamma$ est une racine sph\'erique, alors l'\'el\'ement $s_\gamma \in W^Q$ correspondant v\'erifie :
\begin{liste}
\item[i)] $s_\gamma = s_\alpha s_\beta$ pour certaines racines positives $\alpha , \beta$ telles que $\cg \alpha,\beta^\ch \cd = 0$ et $\alpha + \beta \in \Z_{>0} \gamma$ ;
\item[ii)] $s_\gamma \rho \in  \rho + \Z\gamma$.
\end{liste}
\end{lem}

{\bf Remarque :} L'assertion ii) entra\^ine  $\cg \rho , \alpha^\ch \cd = \cg \rho , \beta^\ch \cd $ car les racines $\alpha$ et $\beta$ ne sont pas proportionnelles et :
$$\rho - s_\gamma \rho  = \rho - s_\alpha s_\beta \rho  =   \cg \rho , \alpha^\ch \cd \alpha + \cg \rho , \beta^\ch \cd \beta \in  \Z (\alpha + \beta) \p$$
 Si $\lambda \in \pic(X)$ alors  on a aussi : $\cg \lambda , \alpha^\ch \cd = \cg \lambda , \beta^\ch \cd $ \`a cause de (\ref{eq:sgl}) (\cf la remarque p. \pageref{eq:sgl}).  

\begin{dem}
Soit $\gamma$ une racine sph\'erique de $X$. 

D'abord, on se ramène au cas o\`u $X$ est de rang $1$ : 

On pose $$X_1:= \inter_{j \not = \gamma} D_j$$
(l'intersection transverse de diviseurs limitrophes $D_j \not= D_\gamma$). Remarquons que $X_1$ est encore magnifique de rang minimal car ses points fixes pour $T$ sont dans l'orbite ferm\'ee $F$ (\cf \cite[pro. 2.3]{Nico}). Au lieu de raisonner dans la vari\'et\'e $X$, on peut raisonner dans la vari\'et\'e $X_1$ qui est de rang $1$ et dont $\gamma$ est la racine sph\'erique. On suppose donc pour la suite de la d\'emonstration que $X$ est de rang $1$.

\begin{center}
{\bf 1ère \'etape : quelques cas particuliers}
\end{center}

Le lemme est v\'erifi\'e pour les exemples de vari\'et\'es magnifiques du paragraphe \ref{ssec:exrun} :

--- l'espace projectif $\PP^{2n-1}$ et la quadrique $\qu_{2n-1}$ munis de l'action naturelle de $Spin_{2n}(\kk)$ ;

--- l'espace projectif $\PP^7$ et la quadrique $\qu_7$ munis de l'action naturelle de $Spin_7(\kk)$.

En effet, voici pour chacun de ces cas-l\`a, dans le tableau suivant (qui complète le tableau du paragraphe \ref{ssec:exrun}),  $s_\gamma \rho - \rho$ et deux racines positives $\alpha,\beta$ qui satisfont le i) du lemme :

$$
\begin{array}{|c|c|c|c|c|}
\hline
\mbox{ diagrammes } &\multicolumn{2}{|c|}{\scriptstyle \dnn \;\;\;\; (n\ge 3) \mbox{ ou } \axaa} & \multicolumn{2}{|c|}{\scriptstyle \spii}\\
\hline
G & \multicolumn{2}{|c|}{Spin_{2n}(\kk)} & \multicolumn{2}{|c|}{Spin_7(\kk)}\\
\hline
X & \qu_{2n-1} & \PP^{2n-1} & \qu_7 & \PP^7 \\
\hline

\alpha \atop
\beta&\multicolumn{2}{|c|}{{\alpha_1 +... + \alpha_{n-1} \atop  \alpha_1+...+ \alpha_{n-2} + \alpha_n }\mathrm{ou} {\alpha \atop \alpha'} \mathrm{si\;} n=2} &\multicolumn{2}{|c|}{\alpha_1+\alpha_2+2\alpha_3 \atop  \alpha_2+\alpha_3}\\
\hline
\gamma& \scriptstyle \frac{\alpha+\beta}{2}& \scriptstyle \alpha+\beta & \scriptstyle \frac{\alpha+\beta}{2}&\scriptstyle \alpha+\beta\\
\hline
s_\gamma&\multicolumn{4}{|c|}{s_\alpha s_\beta}\\
\hline
s_\gamma \rho - \rho&(2-2n) \gamma &(1-n)\gamma &(2-2n) \gamma&(1-n)\gamma\\
\hline
\end{array}
$$

Maintenant, on va r\'eduire la d\'emonstration du lemme \`a ces cas particuliers.
\begin{center}
{\bf 2ème \'etape : induction parabolique}
\end{center}

Rappelons une d\'efinition de \cite[\S 3.4]{L01}

\begin{defin}\label{defi:indupara}
Soit $P$ un sous-groupe parabolique de $G$ de radical $R(P)$. Si $\adh{X}$ est une $P/R(P)-$vari\'et\'e magnifique, alors on note : $G\croi^P \adh{X}$ la $G-$vari\'et\'e alg\'ebrique obtenue comme quotient de $G \croi \adh{X}$ sous l'action de $P$ donn\'ee par :
$$\qq p \in P, \qq g \in G, \qq x \in \adh{X}, p.(g,x) := (gp\inv,p.x) \p$$
La vari\'et\'e $G\croi^P \adh{X}$ est une $G-$vari\'et\'e magnifique et on dit qu'elle est obtenue par induction parabolique de $\adh{X}$ \`a travers $P$.
\end{defin}

{\bf Remarques :} 1) En ce qui concerne le rang, on a l'\'equivalence :

$$G\croi^P \adh{X}\mbox{  est de rang $r$ (resp. de rang minimal)}$$
$$\equi \adh{X} \mbox{ est de rang $r$ (resp. de rang minimal).}$$

2) En fait $\adh{X}$ est une sous-vari\'et\'e ferm\'ee de $G\croi^P \adh{X}$ car pour le morphisme 
$$\theta :  G\croi^P \adh{X} \to G/P \;\;,\;\; (g,x) \mod P \donne g \mod P$$
$\adh{X} = \theta\inv(P/P)$.

Nous allons voir que l'induction parabolique conserve les propri\'et\'es $i), ii)$. 

Consid\'erons donc une vari\'et\'e magnifique $\adh{X}$ de rang minimal $1$ et supposons que $X =G \croi^{Q'} \adh{X}$ pour un certain parabolique $Q'$ de $G$.  D'abord, quitte \`a conjuguer $Q'$ par un \'el\'ement de $G$, on peut exiger que $Q \sub Q'$. En effet, l'existence d'un morphisme $G-$\'equivariant $X = G\croi^{Q'} \adh{X} \to G/{Q'}$ montre que $G/Q'$ contient un point fixe pour $Q$.

On va montrer que si $\adh{X}$ v\'erifie le lemme \ref{lem:sg}, alors $X$ aussi. Fixons pour cela les notations suivantes concernant $\adh{X}$ : soit $R(Q')$ le radical de $Q'$ et 

$$\adh{G} := Q' / R(Q') , \; \adh{T} : = T/R(Q') \cap T , \; \adh{B} : = B / R(Q') \cap B, $$
$$\adh{B^-}:= B^- / R(Q') \cap B^-, \; \adh{Q} := Q/R(Q')\cap Q , \; \adh{W}:= W_{Q'} = N_{Q'}(T)/T \p$$
Dans $\adh{G}$, $\adh{B}$ est un sous-groupe de Borel et $\adh{T}$ un tore maximal. Notons $\adh{\Phi}$ le système de racines correspondant avec ses racines positives $\adh{\Phi}^+$ et ses racines simples $\adh{\Delta}$. On a alors les inclusions :
$$\adh{\Phi} \sub \Phi , \adh{\Phi}^+ \sub \Phi^+, \adh{\Delta} \sub \Delta \p$$

{\bf Remarque :} le groupe de Weyl de $(\adh{G}, \adh{T})$ est $\adh{W}$ qui est un sous-groupe de $W$ et on a aussi : $\adh{W}^{\adh{Q}} \sub W^Q$.

Enfin, notons $\adh{\gamma}$ l'unique racine sph\'erique de $\adh{X}$ (il n'y en a qu'une car $\adh{X}$ est de rang $1$). 

En fait $\adh{\gamma} = \gamma$, la racine sph\'erique de $X$.  
En effet, si on note $\adh{F}$ l'unique $\adh{G}-$orbite ferm\'ee de $\adh{X}$ et $F$ l'unique $G-$orbite ferm\'ee de $X$, alors d'une part $\adh{\gamma}$ est le poids de $\adh{T}$ dans $T_{\adh{\zzz}} \adh{X} /T_{\adh{\zzz}} \adh{F}$ (o\`u $\adh{\zzz}$ est l'unique point fixe de $\adh{B^-}$ dans $\adh{X}$) et d'autre part, $\gamma$ est le poids de $T$ dans $T_{\zzz} X / T_{\zzz} F$ (o\`u $\zzz$ est l'unique point fixe de $B^-$ dans $X$). Mais puisque l'on a suppos\'e que $X = G\croi^{Q'}\adh{X}$, on a $ \zzz = {1}\croi^{Q'} \adh{\zzz} $ et $F = G \croi^{Q'}\adh{F}$ et donc l'isomorphisme de $T-$modules :
$$ T_{\adh{\zzz}} \adh{X} /T_{\adh{\zzz}} \adh{F} \iso T_{\zzz} X /T_{\zzz} F  \p$$  

Maintenant, vu que $\adh{\gamma} = \gamma$, les courbes irr\'eductibles associ\'ees \`a $\gamma$ dans $X$ et \`a $\adh{\gamma}$ dans $\adh{X}$ co\"{i}ncident. On a par cons\'equent : $$s_\gamma = s_{\adh{\gamma}} \in \adh{W}^{\adh{Q}} \sub W^Q \p$$ 

Puisque $\adh{\Phi}^+ \sub \Phi^+$, si $\adh{X}$  v\'erifie le point $i)$ du lemme \ref{lem:sg}, alors $X$ aussi.

Pour le point $ii)$, c'est moins imm\'ediat.

On va utiliser les poids fondamentaux $(\omega_\delta)_{\delta \in \Delta}$.

Comme $\adh{\Delta} \sub \Delta$, on peut \'ecrire :
$$\rho = \sum_{\delta \in \Delta} \omega_\delta = \sum_{\delta \in \adh{\Delta}} \omega_\delta + \sum_{\delta \in \Delta \moins \adh{\Delta}}\omega_\delta $$
$$= \adh{\rho} + \pi$$
o\`u $\adh{\rho}$ est la demi-somme des racines positives de $\adh{\Phi}$ et o\`u $\ma \pi :=\sum_{\delta \in \Delta \moins \adh{\Delta}}\omega_\delta $.

On a alors :

$$s_\gamma (\rho) - \rho = s_{\adh{\gamma}} (\rho) - \rho $$
$$ =  s_{\adh{\gamma}} (\adh{\rho}) - \adh{\rho} + s_{\adh{\gamma}} (\pi) - \pi \p$$
Mais pour chaque $\alpha \in \adh{\Delta}$ :
$$s_{\alpha} \pi = \pi - \cg \pi , \alpha^\ch \cd \alpha$$
$$ = \pi - \sum_{\delta \in \Delta \moins \adh{\Delta}} \cg \omega_{\delta} , \alpha^\ch \cd \alpha$$
$$= \pi \p$$
Donc, pour tout $w \in \adh{W}$, $w \pi = \pi$ et en particulier :
$$s_{\adh{\gamma}} \pi = \pi$$
car $s_{\adh{\gamma}} \in \adh{W}$.

Ainsi, on trouve que :
$$s_\gamma \rho - \rho = s_{\adh{\gamma}} \adh{\rho} - \adh{\rho} \p$$

Finalement, si $\adh{X}$ v\'erifie $ii)$, alors $X$ aussi.

Pour conclure cette d\'emonstration, on utilise que toute vari\'et\'e magnifique de rang $1$ est une induction parabolique
$$G \croi^P \adh{X}$$ 
o\`u $\adh{X}$ est une des vari\'et\'es magnifiques \'etudi\'ees en exemple ci-dessus :
$$\PP^{2n-1} \ou \qu_{2n-1}$$ munie de l'action de $ Spin_{2n}(\kk)$  ou encore $$ \PP^7 \ou \qu_7$$ munie de l'action de $ Spin_7(\kk)$.
 En effet, ce sont les seules vari\'et\'es magnifiques de rang minimal $1$ qui apparaissent dans \cite[th. A]{Nico} (\cf aussi le tableau \cite[table 1, p. 381]{Was} et  \cite[lemme 2.2, d\'ef. 2.3, p. 379]{Was}).
\end{dem}
\section{Retour sur le lemme clef}\label{sec:lemcle}

Rappelons qu'est fix\'e un sous-groupe \`a un paramètre $\zeta : \kk^* \to T$ dominant et $X-$r\'egulier de sorte que :
$$X^\zeta
 = X^T$$
et les cellules positives $X^+(x)$, $x \in X^T$, sont toutes $B-$stables.

\subsection{Rappel des hypothèses}

La situation qui nous int\'eresse (\cf la condition (\ref{eq:cond1}) page \pageref{eq:cond1}) est celle de deux $B-$orbites $\call{B}$ et $\call{B}'$ dans $X$ dans la position relative suivante :
\begin{equation}\label{eq:posrel}
\call{B}' \sub \adh{\call{B}} \et \codim_X\call{B'}  = \codim_X\call{B} +1 = i+1 \p
\end{equation}
et on supppose de plus que le morphisme :
\begin{equation}\label{eq:hypdimu}
d^i_{\call{B},\call{B}'}(\mu) : H^i_{\call{B}}(\call{L}_\lambda)_{(\mu)} \to H^{i+1}_{\call{B}'}(\call{L}_\lambda)_{(\mu)} 
\end{equation}
est non nul.

Notre objectif est de d\'emontrer le lemme clef \ref{lem:memecellule} \cad que $\call{B}$ et $\call{B}'$ sont deux $B-$orbites de la m\^eme cellule.

Pour cela, on note $\call{O}$, $\call{O}'$ les $G-$orbites de $\call{B}$, $\call{B}'$ et $w,w'$ les \'el\'ements de $W^Q$ tels que :
$$\call{B} = X^+_w \cap \call{O} \et \call{B}' = X^+_{w'} \cap \call{O}' \p$$

Rappelons que $X^+_w$ et $X^+_{w'}$ sont les cellules centr\'ees en les points $T-$fixes : $w\zzz$ et $w'\zzz$.

Il s'agit donc de d\'emontrer que $w = w'$ si les hypothèses (\ref{eq:posrel}) et (\ref{eq:hypdimu}) sont v\'erifi\'ees.

\subsection{D\'emonstration du lemme clef}\label{subsec:demlemcle}

On va proc\'eder en plusieurs \'etapes.

Comme $\call{B}' \sub \adh{\call{B}}$, on a forc\'ement : $w' \in \adh{X^+_w}$. Or :

\begin{pro}\label{pro:longfixes}

Soient $w,w' \in W^Q$. 

Si  $w'\zzz \in \adh{X^+_w}$, alors :
$$l(w') \ge l(w) \p$$
\end{pro}

\begin{rem} 
Cela est bien connu dans le cas des vari\'et\'es de drapeaux.
\end{rem}

\begin{dem}

D'apr\`es \cite[th. 1.4 (ii)]{B:inf}, $\adh{X^+_w}$ intersecte $F$, l'orbite ferm\'ee de $X$,  proprement dans $G.\adh{X^+_w}$ \cad : toutes les composantes irr\'eductibles de $$\adh{X^+_w} \cap F$$ sont de dimension $\ma \dim \adh{X^+_w} + \dim F - \dim G.\adh{X^+_w}$.

Or, $\adh{X^+_w \cap F}$ est une de ces composantes (car $X^+_w \cap F$ est un ouvert de $\adh{X^+_w} \cap F$). Donc :
$$\dim  \adh{X^+_w} \cap F = \dim \adh{X^+_w \cap F} = \dim \adh{B w\zzz} \p$$

D'un autre c\^ot\'e :

$$w'\zzz \in \adh{X^+_w} \impliq Bw'\zzz \sub \adh{X^+_w} \cap F$$
d'o\`u : $ \dim \adh{Bw'\zzz} \le \dim \adh{Bw\zzz}$.

\end{dem}

En utilisant que $\call{B}' \sub \adh{\call{B}}$ et aussi que $\codim_X(\call{B}') = \codim_X\call{B} + 1$ on trouve :

\begin{pro}\label{pro:alt}
On a l'alternative suivante :
\begin{equation}\label{eq:ta1}
 l(w') = l(w) + 1 \et \dim \call{O}' = \dim \call{O}
\end{equation}

ou bien

\begin{equation}\label{eq:ta2}
l(w') = l(w) \et \dim \call{O}'+1 = \dim \call{O} \p
\end{equation}
\end{pro}

\begin{dem}

Puisque $\call{B} = X^+_w \cap \call{O}$ et comme $X^+_w$ intersecte proprement $F$ dans $\adh{\call{O}}$ (\cite[th. 1.4 (ii)]{B:inf}), on a les \'egalit\'es :
$$\dim \call{B} = \dim \call{O} - \codim_{\call{O}}(\call{B})$$
$$= \dim \call{O} - \codim_F(X^+_w \cap F)$$
$$= \dim \call{O} - \codim_F(Bw\zzz) $$
$$= \dim \call{O} - l(w) \;\;;$$
et de m\^eme : 
$$\dim \call{B}' = \dim \call{O}'- l(w') \;\;;$$
d'o\`u :
$$\dim \call{O}- \dim \call{O}' + l(w') - l(w) = 1 \p$$

Or, d'une part l'inclusion : 
$$\call{O}' = G. \call{B}' \sub G.\adh{\call{B}} = \adh{\call{O}}$$
entra\^ine que :
$$\dim\call{O} - \dim \call{O}' \ge 0 \;\;;$$
et d'autre part, d'après la proposition \ref{pro:longfixes} :
$$l(w') - l(w) \ge 0\p$$
\end{dem}

Pour les deux dernières propositions, on a seulement utilis\'e l'hypothèse :
$$
\call{B}' \sub \adh{\call{B}} \et \codim_X\call{B'}  = \codim_X\call{B} +1 \p$$
Continuons : puisque $w'\zzz \in \adh{X^+_w}$,  d'apr\`es le lemme \ref{lem:chemin}, il existe, dans la vari\'et\'e  $X$, une cha\^ine de courbes, $c_i$, $T-$invariantes, irr\'eductibles (et ferm\'ees) qui relient $w\zzz$ \`a $w'\zzz$ :
\begin{equation}\label{eq:chemin}
w\zzz= x_0 \sta{{c_0}}{\to} x_1 \sta{c_1}{\to} ... x_N \sta{c_N}{\to} x_{N+1} =w'\zzz \p
\end{equation} 

Cette notation signifie que les courbes $c_i$ \og vont toutes dans le m\^eme sens \fg\ \ie :
$$\qq i ,\; c_i(0) = x_i \et c_i(\infi) = x_{i+1}$$
(o\`u l'on rappelle que pour un point quelconque $y \in c_i \moins c_i^T$, 
$$\ma c_i(0) := \limz \zeta(a).y \et c_i(\infi) := \limi \zeta(a).y $$
pour le sous-groupe \`a un paramètre $\zeta$ qui a servi \`a d\'efinir la d\'ecomposition cellulaire $\ma X = \Dij_{w \in W^Q} X^+_w$).

\vskip .3cm

Comme tout point fixe $x \in X^T$ \'ecrit de manière unique $x =\sigma \zzz$ pour un certain $\sigma \in W^Q$, on notera dor\'enavant :
$$l(x)\index{$l(x)$} : = \codim_F (B.x) = l(\sigma) \p$$ 

On aura besoin du r\'esultat suivant :

\begin{pro}\label{pro:plusquedeux}
Soit $X$ une vari\'et\'e magnifique de rang minimal. 

Alors, pour tout $\lambda \in \pic(X)$ tel que $\lambda +\rho$ est r\'egulier et pour toute racine sph\'erique $\gamma$ de $X$, on a :
\begin{equation}\label{eq:plusquedeux}
|l(w_\lambda) - l(w_\lambda s_\gamma)| = |l(\lambda) - l(s_\gamma * \lambda) | \ge 2
\end{equation}
\end{pro}

\begin{dem}
La première \'egalit\'e a lieu par d\'efinition de $l(\lambda)$, on va d\'emontrer la minoration.

D'après le lemme \ref{lem:sg} (p. \pageref{lem:sg}), il existe $2$ racines positives et orthogonales $\alpha$ et $\beta$ telles que :
$$\alpha + \beta \in \Z \gamma \et s_\gamma = s_\alpha s_\beta \p$$

Or on sait que, pour tout $\lambda \in \pic(X)$, $s_\gamma \lambda - \lambda \in \Z \gamma$ et que $s_\gamma \rho - \rho \in \Z \gamma$ (\cf (\ref{eq:sgl}) p. \pageref{eq:sgl} et le lemme \ref{lem:sg}, ii)). 

Donc :
$$s_\gamma * \lambda - \lambda = s_\gamma (\lambda + \rho) - \lambda - \rho$$
$$= -\cg \lambda +\rho, \alpha^\ch \cd \alpha - \cg \lambda+ \rho , \beta^\ch\cd \beta \in \Z \gamma \p$$

Mais, on a forc\'ement :
\begin{equation}
\cg \lambda +\rho, \alpha^\ch \cd = \cg \lambda+ \rho , \beta^\ch\cd
\end{equation}
(\cf la remarque qui suit le lemme \ref{lem:sg}).

En outre, puisque $\lambda + \rho$ est un caractère r\'egulier, $ \cg \lambda + \rho, \alpha^\ch \cd  > 0$ ou $ \cg \lambda + \rho, \alpha^\ch \cd  < 0$.

Mais alors, le premier cas entra\^ine :
 $$l(s_\gamma * \lambda) = l((s_\alpha s_\beta)* \lambda) > l(s_\beta * \lambda) > l(\lambda)$$
et le deuxième entra\^ine :
$$l(s_\gamma * \lambda) = l((s_\alpha s_\beta)* \lambda) < l(s_\beta * \lambda) < l(\lambda) \p$$
On conclut, dans tous les cas, que :
$$|l(s_\gamma * \lambda) - l(\lambda) | \ge 2 \p$$ 

\end{dem}

\vskip .3cm

Maintenant, on va utiliser l'hypothèse :

$$H^i_{\call{B}}(\call{L}_\lambda)_{(\mu)} \et H^{i+1}_{\call{B}'}(\call{L}_\lambda)_{(\mu)} \not=0 \p$$

On raisonne par l'absurde : on suppose que $w \not=w'$ \ie : la cha\^ine (\ref{eq:chemin}) est non triviale ($N \ge 0$).

Il ne reste plus que trois \'etapes :

\primo dans la cha\^ine (\ref{eq:chemin}), la courbe $c_0$ est contenue dans l'orbite ferm\'ee $F$ ;

\secundo dans la cha\^ine (\ref{eq:chemin}), $N=0$ ;

\tertio Contradiction !
\vskip .5cm

{\bf 1\`ere \'etape :} Pour la suite de points fixes de la cha\^ine de courbes (\ref{eq:chemin}), on a gr\^ace aux propositions  \ref{pro:longfixes} et \ref{pro:alt} les in\'egalit\'es :
\begin{equation}\label{eq:cheminin} 
l(w) = l(x_0) \le l(x_1) \le ... \le l(x_N) \le l(x_{N+1}) = l(w') \le l(w) +1\p
\end{equation}

De (\ref{eq:cheminin}), il d\'ecoule alors que :
\begin{equation}\label{eq:lxx}
0 \le l(x_1) - l(x_0) \le 1 \p
\end{equation}

Si la courbe $c_0$, qui relie $x_0$ \`a $x_1$, n'est pas contenue dans l'orbite ferm\'ee $F$, alors, d'après le paragraphe \ref{subsec:cgamma}, il existe une racine sph\'erique $\gamma_0$ telle que :
$$c_0 = wC_{\gamma_0}$$
($C_{\gamma_0}$ est la courbe irr\'eductible et $T-$invariante associ\'ee \`a $\gamma_0$). Donc : 
\begin{equation}\label{eq:xzxu}
x_0 = w\zzz \et x_1 = ws_{\gamma_0} \zzz \p
\end{equation}

Or, ce $w$ tel que $x_0 = w\zzz$ n'est pas un \'el\'ement quelconque de $W^Q$.

En effet, pour que la composante 
$$H^i_{\call{B}}(\call{L}_\lambda)_{(\mu)}$$
ne soit pas nulle (\ie pour que le $\goth g-$module simple $L(\mu)$ ait une multiplicit\'e non nulle dans le $\goth g-$module $H^i_{\call{B}}(\call{L}_\lambda)$), il est n\'ecessaire, selon le lemme \ref{lem:multi0ou1}, que : 
$$\mu + \rho = w(\lambda + \rho + \sum_\gamma n_\gamma \gamma)$$
o\`u $\ma \sum_\gamma n_\gamma \gamma$ est une combinaison lin\'eaire \`a coefficients entiers de racines sph\'eriques.

En particulier, $w$ est de la forme $w= w_{\lambda'}$ pour un certain $\lambda' = \lambda +\sum_\gamma n_\gamma \gamma \in \pic(X)$ tel que $\lambda' + \rho$ est r\'egulier (rappelons que pour un tel $\lambda'$, $w_{\lambda'}$ est le seul \'el\'ement de $W$ pour lequel $w_{\lambda'}(\lambda' + \rho)$ est dominant et r\'egulier).

On en d\'eduit que :
$ws_{\gamma_0} = w_{\lambda'} s_{\gamma_0} \in W^Q$ : en effet, si $\alpha \in \Phi^+$ est aussi une racine de $Q$, alors :
$$w_{\lambda'} s_{\gamma_0} (\alpha) > 0 \equi (\mu + \rho,s_{\gamma_0} (\alpha) ) > 0$$
$$\equi \left( \,w_{\lambda'}(\lambda' + \rho)\;,\; w_{\lambda'} s_{\gamma_0} (\alpha) \,\right) >0$$
$$\equi \left(\, \lambda' + \rho \; , \; s_{{\gamma_0}}(\alpha) \, \right) >0$$
$$\equi \left(\, s_{{\gamma_0}} \lambda' \;,\;\alpha \, \right) + \left(\, \rho \;,\; s_{\gamma_0}(\alpha) \,\right)  > 0 \p$$

Mais d'une part, comme $s_{\gamma_0} \lambda' \in \pic(X)$ (\cf (\ref{eq:sgl}) dans la remarque p. \pageref{eq:sgl}),  $s_{\gamma_0} \lambda$ est un caractère de $Q$ (et pas seulement de $T$) et donc : $(s_{\gamma_0} \lambda' , \alpha) = 0$. 

D'autre part,  par d\'efinition :
$$s_{\gamma_0} \in W^Q \impliq s_{\gamma_0}(\alpha) >0 \impliq(\rho, s_{\gamma_0}(\alpha) )> 0 \p$$ 

Ainsi, on a montr\'e que $ws_{\gamma_0}(\alpha) >0$ pour toute racine positive $\alpha$ qui est aussi racine de $Q$, \ie : 
\begin{equation}\label{eq:wq}
ws_{\gamma_0} \in W^Q \p
\end{equation}

Il r\'esulte alors de (\ref{eq:xzxu}) et de (\ref{eq:wq}) que :$$|l(x_0) - l(x_1) | = |l(w) - l(ws_{\gamma_0})|$$
$$ = |l(w_{\lambda'}) - l(w_{\lambda'} s_{\gamma_0})| $$
$$ = |l(\lambda') - l(s_{\gamma_0} * \lambda') |\ge 2$$
comme le montre la proposition \ref{pro:plusquedeux}, car $\lambda'+ \rho$ est r\'egulier.
Cela contredit (\ref{eq:lxx}) et en cons\'equence, la courbe $c_0$ est contenue dans l'orbite ferm\'ee $F$.

\saut

{\bf 2\`eme \'etape :} $N=0$.

On d\'eduit de la 1ère \'etape que :
$$\dim B.x_0 \not= \dim B.x_1 \mbox{ \ie } l(x_0) \not= l(x_1)$$
\noindent car maintenant, $x_0$ et $x_1$ sont les deux point fixes d'une courbe $T-$invariante dans $F$ qui est une vari\'et\'e de drapeaux et $B.x_1 \sub \adh{B.x_0}$. 

En cons\'equence :
$$l(x_1) - l(x_0) = 1 \p$$

Pour les m\^emes raisons, la courbe $c_N$ est incluse dans $F$ et :
$$l(x_N) - l(x_{N+1}) = 1 \p$$

Mais d'apr\`es les in\'egalit\'es (\ref{eq:cheminin}), cela n'est possible que si $N=0$.

\saut
{\bf 3\`eme \'etape :} Forc\'ement $w=w'$.

D'apr\`es la 2\`eme \'etape, tout se passe dans la vari\'et\'e de drapeaux $F \iso G/Q$ o\`u $\zzz$ est le point $Q/Q$ et :
$$x_0= w\zzz \;,\; x_1 =w'\zzz \;,\; l(w') = l(w)+1$$
et surtout, les points $w\zzz$ et $w'\zzz$ sont reli\'es par une courbe irr\'eductible et $T-$invariante de $G/Q$. Comme dans $G/Q$, les cellules de Bialynicki-Birula et les $B-$orbites co\"incident, on en d\'eduit que $Bw'\zzz$ est une $B-$orbite de codimension $1$ dans $\adh{Bw\zzz}$. 
Il existe par cons\'equent une racine $\alpha \in \Phi$ telle que :\begin{equation}\label{eq:ws}
w' = ws_\alpha \p  
\end{equation}

Or, puisque la composante $\ma H^i_{\call{B}}(\call{L}_\lambda)_{(\mu)}$ est non nulle, d'apr\`es le lemme \ref{lem:multi0ou1}, il existe un certain poids $\ma a \in \sum_{\gamma \in \Sigma_X} \Z \gamma$, combinaison \`a coefficients entiers de racines sph\'eriques, tel que :
\begin{equation}\label{eq:a}
w\inv(\mu + \rho) = \lambda + \rho +a
\end{equation}
De m\^eme, $w'$ v\'erifie :
\begin{equation}\label{eq:a'}
{w'}\inv (\mu +\rho) = \lambda + \rho + a'
\end{equation}
pour un $\ma a' \in \sum_{\gamma \in \Sigma_X} \Z \gamma$.  

En faisant la diff\'erence (\ref{eq:a'}) - (\ref{eq:a}) et en utilisant (\ref{eq:ws}), on trouve :
$$s_\alpha w\inv (\mu + \rho) - w\inv(\mu + \rho) = a' -a \;,$$
$$ \ie \cg \mu + \rho , (w\alpha)^\ch \cd \alpha = a' - a \;, \mbox{ d'o\`u :}$$
\begin{equation}
 \alpha \in \sum_{\gamma \in \Sigma_X} \Q \gamma
\end{equation}
car, le caract\`ere $\mu + \rho$ \'etant dominant r\'egulier, $(\mu + \rho , w \alpha) \not= 0$.

Mais cela est {\it impossible} ! car comme $X$ est magnifique de rang minimal, on a : 
$$\Phi \inter  \left(\sum_{\gamma \in \Sigma_X} \Q \gamma  \right)= \vide$$
(\cf le lemme \ref{lem:racnonsph} de la section \ref{sec:rac-sph}). \qed

\begin{center}
***
\end{center}



{\it Je remercie tout particulièrement Michel Brion pour ses remarques et commentaires très clairs ansi que  Nicolas Ressayre pour des discussions très utiles.}

\backmatter

\printindex


\clearpage


\end{document}